\documentclass[]{amsart}
 \usepackage[]{amsmath, amsthm, amsfonts}
 \usepackage[all]{xy}
 \usepackage{hyperref}

 \newcommand {\N} {{\mathbb N}}
 \newcommand {\C} {{\mathbb C}}
 \newcommand {\R} {{\mathbb R}}
 \newcommand {\Z} {{\mathbb Z}}
 \newcommand {\Q} {{\mathbb Q}}
 \newcommand {\PP} {{\mathbb P}}
 
 \newcommand {\F} {{\mathcal F}}

 \newcommand {\E} {{\mathcal E}}
 \newcommand {\dt} {{\bullet}}

 \newcommand {\Hbm} {H}

 \newtheorem{thm}[subsection]{Theorem}
 \newtheorem{cor}[subsection]{Corollary}
 \newtheorem{lemma}[subsection]{Lemma}
 \newtheorem{prop}[subsection]{Proposition}
 
 \newtheorem{remark}[subsection]{Remark}

\begin{document}

\title{
Hodge cycles on some moduli spaces }

\author{ Donu Arapura}
\thanks{Author partially supported by the NSF}
\address{Department of Mathematics\\
Purdue University\\
West Lafayette, IN 47907\\
U.S.A.}

\maketitle

The primary goal of this paper is to verify  the Hodge and
generalized Hodge conjectures  for  certain  moduli spaces of sheaves 
over curves and surfaces. When $X$ is a smooth projective curve, let
$U_X(n,d)$ be the moduli space of stable vector bundles of
rank $n$ and degree $d$. If $n$ and $d$ are coprime, this space is
known to be smooth and projective. del Ba\~no has shown
that the Hodge  conjecture
holds for $U_X(n,d)$ if it holds for all powers of the Jacobian
$J(X)$. We reprove this along with an analogous statement for the
generalized Hodge conjecture.
From this, it
is easy to deduce a refinement of  result of Biswas and Narasimhan \cite{BN}
that the generalized Hodge conjecture is valid for 
$U_X(n,d)$ when $X$ is very general
in moduli. We have some extensions of these results when
$X$ is a smooth projective surface.
In particular, we show that the Hodge conjecture holds for the moduli
space of semistable torsion free sheaves over an abelian or Kummer surface.
For arbitrary surfaces,  these spaces are difficult to analyze, so
our attention is devoted to the rank one case,  i.e. the Hilbert
scheme of points.
We   show that the Hodge (respectively generalized Hodge) conjecture holds
for the Hilbert scheme of points on $X$ if it holds for all powers of
$X$. For the Hodge conjecture, this last result can be deduced from 
some work of de Cataldo and  Milgiorini \cite{cm1,cm2}.

The key reductions in the proofs of these theorems
are based on some simple criteria established in the first section.
Here we show that  the Hodge conjecture
holds for a  smooth projective variety if can be
dominated by, stratified by, or fibered over (with
suitable fibers)  varieties where the conjecture
holds. The second case forces us to deal with the
Hodge conjecture for quasiprojective varieties;
we use Jannsen's formulation of it. The second section
contains a few applications of these ideas apart from
the main theorems. For example that the Hodge conjecture
holds for a smooth projective variety with $\C^*$-action if it
holds for the fixed point locus. The third section collects
some results pertaining to the Hodge conjecture
for powers of varieties. These along with the previous
lemmas yield the main theorems.
In the final section, we consider some arithmetic analogues
of these results. We give criteria for
the validity of Tate's conjecture for some of the spaces
considered above. Also we show, in accordance with a
conjecture (or ``espoir'') of Deligne, that Hodge
cycles are absolute on the moduli space of vector bundles over any
curve, or on the Hilbert scheme of points over any surface of
Kodaira dimension zero.

Except for the last section, all varieties will be
defined over $\C$, and (co)homology groups will 
be with respect to the  usual topology with rational
coefficients.
Following Koll\'ar, we will refer to a point of
the  complement of a countable union of proper 
analytic subvarieties  of an analytic variety 
as a {\em very general} point.
If a variety has an obvious moduli space,
we say that the variety is very general if the
corresponding point moduli is very general.

I would like to thank M. de Cataldo, E. Markman, K. Murty, M. Nori, P. Sastry,
and Y. Zarhin for providing  useful feedback and references. 
Parts of this paper were written  during visits
to the University of Chicago in the fall of 2000 and MSRI in the
April 2002. My thanks to these institutions for their hospitality.

\section{Basic tools}

We first  recall some basic notions from Hodge theory.
The primary reference for mixed Hodge structures is
\cite{deligne-hodge}.
We will refer to a pure rational Hodge structure simply as a
Hodge structure.
Given a smooth projective variety $X$ and  codimension $p$
subvariety $Z\subset X$, it will be convenient to
view the fundamental class $[Z]$ as an element of the weight $0$ Hodge
structure $H^{2p}(X,\Q)(p)$. The $(p)$ indicates the
Tate twist  of the Hodge structure \cite[2.1.13]{deligne-hodge}; this amounts
to shifting the Hodge bigrading by $(-p,-p)$ and modifying
the lattice by a factor of $(2\pi i)^p$.
Let $H^{2p}_{alg}(X)$ denote the $\Q$-span of these classes.
For any weight $0$ Hodge structure $H$, let $H_{hodge}$
denote the the intersection of $H^{00}$ with
the rational lattice. We will write $H_{hodge}^{2p}(X)$ for
$[H^{2p}(X)]_{hodge}$. We always
have an inclusion $H^*_{alg}(X)\subset H_{hodge}^*(X)$;
the Hodge conjecture asserts the converse.

Next, we  recall the generalized Hodge conjecture.
This requires some additional notation.
The level of a Hodge structure $H = \oplus H^{pq}$ is the
maximum of $\{|p-q|\,|\, H^{pq}\not= 0\}$. 
Given a complex subspace $W\subseteq V$ of a rational Hodge structure,
let $W_h$ be 
largest sub Hodge structure of $V$ contained in $W$. Note that any sub 
Hodge structure, such as
$W_h$, is determined by  the underlying rational subspace, and 
we will usually identify it with the subspace.
Given a Hodge structure $V$, we get a filtration
$\F^pV= (F^pV)_h$ by sub Hodge structures that we will call the level 
filtration. If $V$ has weight $m$,
$\F^pV$ is precisely the largest sub Hodge structure
of level at most $|m-2p|$.

Given a smooth projective variety $X$,
the coniveau, or arithmetic filtration, is given by 
\begin{eqnarray*}
N^p H^i(X,\Q) &=& \sum_{codimY\ge p} ker[H^i(X,\Q)\to H^i(X-Y,\Q)]\\
&=& \sum_{codim Y=q\ge p} im[H^{i-2q}(\tilde Y,\Q)(-q)\to H^i(X,\Q)]
\end{eqnarray*}
where $Y$ ranges over closed subvarieties; in the second
expression $\tilde Y\to Y$ are chosen desingularizations and the
maps on cohomology are the Gysin maps. Since the level of $H^{i-2q}(\tilde Y,\Q)(-q)$
is bounded by $i-2p$, we have an inclusion
$$N^p H^i(X,\Q) \subseteq \F^pH^i(X,\Q)$$
Grothendieck's amended version of the generalized Hodge conjecture
 $GHC(H^i(X),p)$ asserts that equality holds \cite{groth2}.
We will say that the generalized Hodge conjecture holds for $X$ if
$GHC(H^i(X),p)$ is true for all $i$ and $p$.
The space $N^pH^{2p}(X)$ is just the subspace generated
algebraic cycles of codimension $p$, while $\F^pH^{2p}(X)=
F^pH^{2p}(X)$ is  the space of Hodge cycles i.e. rational $(p,p)$
classes. Hence, $GHC(H^{2p}(X),p)$ is  the
usual Hodge conjecture. 

It will be convenient to define
$$N^p(H^i(X)(c)) = N^{p+c}H^i(X).$$
The notation is chosen so that the inclusion $N^p\subseteq \F^p$
persists after Tate twisting.

We turn now to the functoriality properties of the coniveau
filtration.

\begin{prop} The filtrations $N^\dt$ and $\F^\dt$ are
preserved by 
\begin{enumerate}
\item pushforwards: if $f:X\to Y$ is a map of smooth projective
  varieties of dimensions $n$ and $m$ respectively,
 then $$f_*(N^pH^i(X))\subseteq N^p(H^{i+2(m-n)}(Y)(m-n))$$

\item pullbacks: if $f$ is as above, then
$$f^*(N^pH^i(Y))\subseteq N^pH^i(X),$$ and 
\item products:
$$N^p(H^i(X))\otimes N^q(H^j(Y))\subseteq N^{p+q}H^{i+j}(X\times Y)$$
\end{enumerate}
\end{prop}

\begin{proof} 
An element $t\in N^pH^i(X)$ lies in the image of a map
$k_*H^{i-2q}(T)(-q)$ where $k:T\to X$ is a morphism
from a smooth  projective variety of dimension $n-q\le n-p$.
Therefore 
$$f_*(t)\in (f\circ k)_*H^{i-2q}(T)(m-n-q)
\subseteq N^p(H^{i+2(m-n)}(Y)(m-n))$$
This proves the first part.

For the third statement. Let $T\to X$ and $S\to Y$ be morphisms from 
smooth projective varieties such that $dim T\le dim X-p$ and
$dim S\le dim Y-q$. Then $dim T\times S\le dim X\times Y -p-q$.
It follows that the image of $T\times S$ lies in
$N^{p+q}(H^{i+j}(X\times Y)$ as expected.

We now turn to the proof of the second part which is the most
involved. To avoid excessive notation, we will suppress Tate
twists.
To begin with, let us assume that $f$ is surjective. Suppose that
the $S\subset Y$ is   an irreducible codimension $q\ge p$ subvariety. 
The preimage $f^{-1}S$ will have codimension 
less than or equal to $q$. By taking general hyperplane sections, we
can find a cycle $Z\subseteq f^{-1}S$ of codimension exactly $q$
surjecting onto $S$. By stratification theory \cite[pp. 33-43]{gm}, we can find
a proper Zariski closed set $Z''\subset Z$ containing the union of singular loci
$Z_{sing}\cup f^{-1}S_{sing}$, such that  the map $f:X-Z''\to Y-f(Z'')$ is
locally trivial along tubular neighborhoods of $Z'=Z-Z''$ and $S'=S-f(Z'')$.
To make the last condition precise, consider the diagram
$$
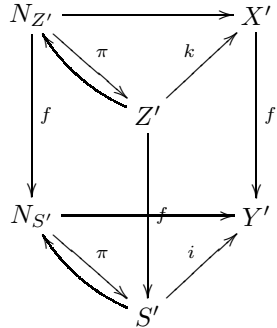
\xymatrix{
 N_{Z'}\ar[rr]\ar[rd]^{\pi}\ar[dd]^{f} &  & X'\ar[dd]^{f} \\
  & Z'\ar[ru]^{k}\ar[dd]^{f}\ar@/^/[lu] &  \\
 N_{S'}\ar[rr]\ar[rd]^{\pi} &  & Y' \\
  & S'\ar[ru]^{i}\ar@/^/[lu] & 
}
$$
where $X' = X - Z''$, $Y'= Y - f(Z'')$, and  $N_{Z'}$
$N_{S'}$ denotes appropriately chosen tubular neighbourhoods of $Z'$
and $S'$ respectively.  The above
condition  is that $N_{Z'}\to f^*N_{S'}$ is a locally trivial map
of locally trivial  fiber bundles (for the classical topology) over
$Z'$.
Fiberwise, we have an open immersion of $2q$ real dimensional oriented
manifolds, and this induces an isomorphism of compactly supported
$2q$ dimensional cohomologies. 
Thus the Thom class $\tau_{Z'}$ of  $N_{Z'}$, which can be viewed element
of relative cohomology
$H^{2q}_{Z'}(N_{Z'} )
 = H^{2q}(N_{Z'}, N_{Z'}-Z')$,
coincides with the pullback of the Thom class $\tau_S'$ on $f^*N_{S'}$. 
The Gysin map $H^*(Z')\to H^{*+2q}(X')$ is given by $\alpha\mapsto
\pi^*\alpha\cup\tau_{Z'}$ extended by $0$
to $X'$. A similar description
holds for $(S',Y')$. It follows that we have a commutative diagram
$$
\xymatrix{
 H^*(S')\ar[r]\ar[d] & H^{*+2q}(Y')\ar[d] \\
 H^*(Z')\ar[r] & H^{*+2q}(X')
}
$$
Therefore, if 
$$\alpha \in ker[H^*(Y)\to H^*(Y-S)= H^*(Y'-S')],$$
 then its image in $H^*(X)$ maps to 
$$ker[H^*(X)\to H^*(X-Z)=H^*(X'-Z')].$$
This implies that $f^*$ preserves $N^p$. The general case is carried
out in the same way after replacing $S$ by $S\cap f(X)$.

\end{proof}  

A correspondence is an algebraic cycle on $Y\times Z$.
Suppose that $T$ is a pure codimension $c$ correspondence which
defines an element $[T]\in H^{2c}(Y\times Z)(c)$.
This induces a morphism $T_*:H^*(Y)\to H^{*+2(c-d)}(Z)(c-d)$
given by $\alpha\mapsto p_{Z*}(p^*_Y(\alpha)\cup [T])$,
where $d =dim \, Y$.

\begin{cor}\label{cor:corr}
 The action of a correspondence
preserves  the above filtrations.
\end{cor}

\begin{lemma}\label{lemma:polar}
 The operation $V\mapsto \F^pV$  induces an
 exact functor from the category of polarizable
 Hodge structures of a fixed weight to itself
\end{lemma}

\begin{proof}
    Let the weight be $m$.
 First note that $\F^p$ is a functor: given a morphism
$f:V\to W$, $f(\F^pV)$ is a sub Hodge structure of
$W$ lying in $F^pW$, since its level  is bounded by $|m-2p|$.
Thus  $f(\F^pV)$ lies in $\F^pW$.
 
Suppose that
$$0\to U\to V\to W\to 0$$
is an exact sequence of polarizable Hodge structures.
By \cite[2.3.5]{deligne-hodge},
$$0\to F^pU\to F^pV\to F^pW\to 0$$
is exact. Certainly, this yields a complex
$$\F^pU\to \F^pV\to \F^pW$$
which will be shown to  be a short exact  sequence.
Injectivity of the first map above is clear. The kernel
of $\F^pV\to \F^pW$ is a Hodge structure lying
in $F^pU$. Thus it must coincide with $\F^pU$.
It remains to check surjectivity of the map $\F^pV\to \F^pW$.
The category of polarizable Hodge structures is semisimple \cite[4.2.3]{deligne-hodge}.
Therefore, there is a splitting $s:W\to V$ for the morphism 
$V\to W$. From the first paragraph we get $s(\F^pW)\subseteq \F^p(V)$.
This finishes the proof.
\end{proof}    

Putting these results together leads to one of the main tools of 
this paper.

\begin{cor}\label{cor:lifthodge} Let $X$ and $Y$ be smooth projective varieties
and suppose that $f:H^{i}(X)\to H^{j}(Y)((j-i)/2)$ is a surjective
morphism of
Hodge structures induced by a correspondence (note that $(j-i)/2$
would be an integer).
Then $GHC(H^i(X),p)$ implies $GHC(H^j(Y), p+(j-i)/2)$.
\end{cor}

\begin{proof}
Let $e= (j-i)/2$, and assume $GHC(H^i(X),p)$.
 The lemma implies that any element $\alpha \in \F^{p+e}H^{j}(Y)(e)$ can be
 lifted to a class $\beta\in \F^pH^{i}(X)=N^pH^i(X)$. Therefore
 $\alpha$ lies in $N^{p+e} $.
\end{proof}

  Next we want to consider the problem of checking the Hodge
 conjecture for a stratified variety. We will reduce it  to the
 conjecture for  strata. 
This forces us to deal with Hodge cycles on arbitrary quasiprojective 
varieties. We recall the formulation of these notions due to Jannsen
\cite{jannsen}.
The appropriate setting for this is homology.
For our purposes the Borel-Moore homology $\Hbm_i(U)$
of  complex algebraic variety $U$ can be taken to be
the  dual of the compactly supported cohomology $H^i_c(U)$.
 $\Hbm_i(U)$ carries a mixed Hodge structure dual to
the one on $H^i_c(U)$ constructed by Deligne \cite{deligne-hodge}.
 From this it follows easily that the weights of $\Hbm_i(U)$
are concentrated in the interval $[-i,0]$.
This mixed Hodge structure is polarizable in the sense
that the associated graded  with respect to the weight
filtration is polarizable.
  We define the space
of Hodge cycles in $\Hbm_{2i}(U)$ to be 
\begin{eqnarray*}
{\Hbm}^{hodge}_{2i}(U) &= &Hom(\Q(0), \Hbm_{2i}(U)(-i)) \cong Hom(\Q(i), W_{-2i}\Hbm_{2i}(U)) \\
&\cong& W_{-2i}\Hbm_{2i}(U)\cap F^{-i}\Hbm_{2i}(U,\C).
\end{eqnarray*}

Given a closed irreducible $i$-dimensional subvariety $V\subset U$,
there is a fundamental class $[V]\in \Hbm^{hodge}_{2i}(U)$.
The space of algebraic cycles $\Hbm^{alg}_{2i}(U)$ is the
span of these classes.
The Hodge conjecture for $U$ in degree $2i$ asserts that $\Hbm^{alg}_{2i}(U)
= \Hbm^{hodge}_{2i}(U)$.
This can be extended to the generalized Hodge conjecture.
We define the niveau filtration by 
$$ N_p\Hbm_i(X,\Q) = \sum_{dim Y\le p} im[\Hbm_i(Y,\Q)\to \Hbm_i (X,\Q)]$$
This is a filtration by submixed Hodge structures.
We will only be concerned with its intersection with  the pure
Hodge structure $W_{-i}\Hbm_i (X,\Q)$. 
Fix a compactification $\bar X$.
Then an easy argument involving weights shows that
$$ N_pW_{-i}\Hbm_i(X,\Q) \cong 
\sum_{dim Y\le p} im[\Hbm_i(\tilde Y,\Q)\to W_{-i}\Hbm_i (\bar X,\Q)]$$
as $\tilde Y$ varies over desingularitations of closed subvarieties.
This together with  Poincar\'e duality:
$$\Hbm_i(\tilde Y,\Q)\cong H^{2dimY-i}(\tilde Y, \Q)(dim Y)$$
gives an estimate on the level which shows
$$N_{p}W_{-i}\Hbm_i(X,\Q) \subseteq
\F^{-p}W_{-i}\Hbm_i(X,\Q)=\F_pW_{-i}\Hbm_i(X,\Q) $$
(where $\F_p$ is defined by the last equality).
Following Jannsen \cite{jannsen} and Lewis \cite{lewis}, we say that
$GHC(H_i(X),p)$ holds if equallity holds above. 
Note that if $X$ is smooth and projective
of dimension $n$, then $GHC(H_i(X), p) $ is equivalent to 
$GHC(H^{2n-i}(X), n-p)$ by Poincar\'e duality. It turns out that the
generalized Hodge conjecture in this setting 
is no stronger than in the usual formulation.

\begin{prop}\label{prop:jannsen}
    If $X'$ is a desingularization of a compactification
 of $X$, then $GHC(H_i(X'),p) $ implies $GHC(H_i(X),p)$.
\end{prop} 

The proof hinges on  natural isomorphism
$$W_{-i}\Hbm_i(X,\Q)\cong \Hbm_{-i}(X',\Q),$$
see \cite[7.9]{jannsen} and \cite[appendix A]{lewis}.

Note that $\Hbm_*$ is functorial in two different ways. It is
covariant for proper maps, and contravariant for open immersions.
These are dual to the proper restriction and extension by zero
maps on $H_c^*$. It follows that the mixed Hodge structure on $\Hbm_*$ 
is compatible with both kinds of maps. It is also easy to see that
the niveau filtration is also compatible with these.

\begin{lemma}\label{lemma:strat}
Let $Z\subset X$ be a closed subset of a projective variety
$X$, and let $U = X-Z$. Then
 $GHC(H_i(U), p)$ and $GHC(H_i(Z),p)$ imply
$GHC(H_i(X),p)$.
\end{lemma}

\begin{proof}
The  exact sequence 
$$\Hbm_{i}(Z)\to \Hbm_{i}(X)\to \Hbm_{i}(U)\to \Hbm_{i-1}(Z)\ldots $$
is compatible with mixed Hodge structures because it is dual
to the long exact sequence for the cohomology of a pair.
Since $H\mapsto W_* H$ is exact on the category of mixed Hodge
structures \cite[2.3.5]{deligne-hodge}, we get an exact sequence
of pure polarizable Hodge structures
$$W_{-i}\Hbm_{i}(Z)\to W_{-i}\Hbm_{i}(X)\to W_{-i}\Hbm_{i}(U).$$
Lemma \ref{lemma:polar} implies that 
$$\F_p W_{-i}\Hbm_{i}(Z)\stackrel{f}{\to} \F_p W_{-i}\Hbm_{i}(X)
\stackrel{g}{\to} \F_pW_{-i}\Hbm_{i}(U)$$
is exact.

Suppose that $\alpha\in \F_p W_{-i}\Hbm_{i}(X)$.
 The hypothesis implies that $g(\alpha)\in N^p W_{-i}\Hbm_{i}(U)$.
Thus $g(\alpha)$ is the  image $\beta\in W_{-i}\Hbm_i(T)$  under
 the map on homology induced by a  proper map $T\to U$,
where $T$ is a finite union of varieties of dimension $\le p$.
Choose a compactification $\bar T $ of $T$. By blowing up, if
 necessary, we can assume that there is a map $\bar T\to X$ extending $T\to U$.
As above, we have a surjection  
$$W_{-i}\Hbm_i(\bar T)\to W_{-i}\Hbm_i(T).$$
Let $\beta'$ be an element of the  space on left lifting $\beta$.
Then $\beta'$ maps to an element $\alpha'\in N_pW_{-i}\Hbm_i(X)$
such that the difference $\alpha-\alpha'$ lies in the 
kernel of $g$, and hence is given by $f(\gamma)$ where
$\gamma\in \F_p W_{-i}\Hbm_{i}(Z)$. By hypothesis
$\F_p W_{-i}\Hbm_{i}(Z)= N_p W_{-i}\Hbm_{i}(Z)$.
Therefore $\alpha = \alpha'+\gamma\in N_p W_{-i}\Hbm_{i}(X)$.

\end{proof}

 By a stratification of an algebraic variety $X$, we will mean
a finite partition $X=\cup X_i$ into locally closed sets, called
strata, such that closure of any stratum is a union of strata.

\begin{cor}\label{cor:strat}
Suppose that $X$ is a projective variety.
If $X$ has a stratification such that each stratum $S$ 
satisfies $GHC(H_i(S), p)$, then  $GHC(H_i(X), p)$
holds.
\end{cor}

\begin{proof} Since a stratum of minimal dimension is closed,
the result follows by induction on the number of strata.
\end{proof}

\begin{lemma}\label{lemma:fiber1}
Let $X = Y\times F$ such 
 that $\Hbm_*(F)$ is spanned by algebraic
cycles. Then  $GHC(H_{i-2j}(Y),p-j)$ for all $j$, such that $i-2j\ge 0$
and $H_{2j}(F,\Q)\not= 0 $,
 implies $GHC(H_i(X),p)$.
\end{lemma}

\begin{proof}
The K\"unneth formula implies that
$$ \Hbm_{k}(X) = \bigoplus_{i+2j=k}\, \Hbm_{i}(Y)\otimes \Hbm_{2j}(F)$$
since by assumption the odd degree homology of $F$ vanishes.
We can choose a  basis of $\Hbm_{2j}(F)$ consisting of fundamental
classes of varieties. These classes are pure of type $(-j,-j)$,
so $\Hbm_{2j}(F)$ is a sum of $\Q(j)$'s.
Since the  K\"unneth decomposition respects mixed Hodge structures, we
have
$$ \F_pW_{-k}\Hbm_{k}(X) = \bigoplus_{i+2j=k}\, \F_{p-j}W_{-i}\Hbm_{i}(Y)\otimes \Hbm_{2j}(F)$$
By assumption, the right hand sum  equals
$$ \bigoplus_{i+2j=k}\, N_{p-j}W_{-i}\Hbm_{i}(Y)\otimes
N_jW_{-j}\Hbm_{2j}(F)\subseteq N_pW_{-k}\Hbm_{k}(X$$
\end{proof}

\begin{lemma}\label{lemma:fiber}
 Let $f:X\to Y$ be  a morphism of smooth varieties
which is a Zariski locally trivial fiber
 bundle with fiber $F$. Suppose that $F$ is  smooth
and  that $\Hbm_*(F)$ is spanned by algebraic
cycles. Then 
 $GHC(H_{i-2j}(Y),p-j)$ for all $j$, such that $i-2j\ge 0$
and $H_{2j}(F,\Q)\not= 0 $, implies $GHC(H_i(X),p)$. 
\end{lemma}

\begin{proof}
The argument is similar to the proof of the previous lemma,
however we will work in cohomology. This is justified by
the fact that our spaces are smooth.
 Let $\beta_l'$ be algebraic cycles in $F$ whose fundamental
classes give a basis for $H^*(F)$.
Choose a  Zariski open set $U\subset Y$ such that $f^{-1}U\to U$
is a product. Let $\beta_l$ be the fundamental classes of
the closures in $X$
of the pullbacks of the $\beta_l'$ to $U$. These classes
define a splitting of the restriction $H^*(X)\to H^*(F)$, so we can
apply the Leray-Hirsch theorem \cite[p. 258]{spanier}
to obtain a decomposition
$$ H^{k}(X) = \bigoplus_{i+2j=k}\, H^{i}(Y)\otimes H^{2j}(F).$$
The rest of the proof proceeds exactly as before.
\end{proof}

\begin{cor}\label{cortolemma:fiber}
 Let $f:X\to Y$ be  a morphism of smooth varieties
which is a (Zariski) locally trivial fiber
 bundle with fiber $F$. Suppose that $F$ is  smooth
and  that $\Hbm_*(F)$ is spanned by algebraic
cycles. Then  the Hodge conjecture
holds for $X$ if it holds for $Y$.
\end{cor}

\begin{remark}\label{remark:cell}
A variety has a cellular
decomposition if it has a stratification such that
the strata are affine spaces.
If $F$  has a cellular decomposition, then it satisfies the
hypothesis of the lemma \cite[19.1.11]{fulton}. 
Flag manifolds and nonsingular projective toric
varieties have cellular decompositions.
\end{remark}

We need a variation on corollary \ref{cor:lifthodge}.

\begin{lemma} Suppose that $f:X\to Y$ is a proper
morphism of varieties such that
$f_*:\Hbm_{k}(X)\to \Hbm_{k}(Y)$ is surjective.
Then $GHC(H_k(X),p)$ implies $GHC(H_k(Y),p)$.
\end{lemma}

\begin{proof}
Assume  $GHC(H_k(X),p)$.
As in the proof of lemma \ref{lemma:strat},
we get a surjection of pure polarizable
Hodge structures $W_{-k}\Hbm_{k}(X)\to W_{-k}\Hbm_{k}(Y)$. Consequently,
lemma \ref{lemma:polar} implies that 
$N_pW_{-k}\Hbm_{k}(X)=\F_pW_{-k}\Hbm_{k}(X)$ surjects onto  $\F_pW_{-k}\Hbm_{k}(Y)$.
So $GHC(H_k(Y),p)$  follows.
\end{proof}

We will give a simple criterion for surjectivity of $f_*$. 
Recall that a locally compact Hausdorff space $X$ is a rational
homology manifold of (real) dimension $n$ if for each $x\in X$,
the rational homology
groups  $H_*(X,X-\{x\})\cong H_*(\R^n,\R^n-\{0\})$.
For such a space  the groups $H_n(X, X-\{x\})$ form a local
system as $x$ varies. The space is orientable if this
local system is constant.
Orientable rational homology manifolds 
satisfy Poincar\'e duality with rational coefficients.
The most interesting examples of  orientable rational homology manifolds for us
are varieties with at worst finite quotient singularities.

\begin{lemma}\label{lemma:rathomology}
 If $f:X\to Y$ is a proper surjective map of quasiprojective
varieties such
that $Y$ is  also a rational homology manifold, then $f_*:\Hbm_*(X)\to
\Hbm_*(Y)$ is surjective.
\end{lemma}

\begin{remark}
The properness is necessary to insure that $f_*$ is defined.
\end{remark}

\begin{proof} 
Suppose that $X$ is smooth and that $dim X = dim Y$. Then Poincar\'e
duality gives a map $f^*:\Hbm_*(Y)\to \Hbm_*(X)$ such that
$f_*f^* = deg(f)$, and this implies surjectivity of $f_*$.

In general, an intersection of a finite number of very ample
divisors will produce a subvariety $X'\subseteq X$ such
that $dim X' = dim X$ and $f|_{X'}$ is surjective. Let $X''\to X$
be a desingularization, and let $g:X''\to Y$ be the natural
map. Then $g_*$ is surjective as above, but it factors through
$f_*$, so we are done.
\end{proof}

\begin{cor}\label{cor:lifthodge2}
With the notation of the lemma, if $X$ satisfies the 
Hodge (respectively generalized Hodge) conjecture then so does $Y$. 
\end{cor}

\section{Some simple examples}

Let us work out some examples to illustrate
the previous techniques.

\begin{lemma}
Let $X$ be a smooth projective variety, and let $V\subset X$ be
a smooth closed subvariety. If the $X$ and $V$ satisfy the Hodge
conjecture, then so does the blow up $\pi:\tilde X\to X$ of
$X$ along $V$.
\end{lemma}

\begin{proof}
Let $U= X-V$ and let $E = \pi^{-1}V$. $E\to V$ is a locally
trivial bundle with projective spaces as fibers. Therefore
$E$ satisfies the Hodge conjecture by lemma \ref{lemma:fiber}.
$U$ satisfies the conjecture by proposition \ref{prop:jannsen}.
Therefore, we are done by lemma \ref{lemma:strat}.
\end{proof}

\begin{cor}\label{cor:dim5}
The Hodge conjecture is birationally invariant in dimensions
up to $5$. In other words, if $X$ and $X'$ are two birational
smooth projective varieties with dimension $\le 5$, then
one of them satisfies the Hodge conjecture if and only
if  the other does.
\end{cor}

\begin{proof}
Suppose that $X$ satisfies the Hodge conjecture.
It is well known that the Hodge conjecture holds in dimensions
up to $3$ (by the Lefschetz $(1,1)$ and
hard Lefschetz theorems). Thus the conjecture holds for
a single blow up of $X$, since the center has dimension $\le 3$.
By iteration, the same goes for any $X''\to X$ given
by a sequence of blow ups. Since any $X'$ birational to $X$
is dominated by such an $X''$, the result follows by
corollary \ref{cor:lifthodge2}. 
\end{proof}

\begin{lemma}[Conte-Murre \cite{conte-murre}] A projective uniruled $4$-fold satisfies the Hodge
conjecture.
\end{lemma}

\begin{proof}
Suppose $X$ is a  projective uniruled $4$-fold. Then any desingularization
is also uniruled. Therefore we can assume that $X$ is smooth by
proposition \ref{prop:jannsen}. There exists a surjective
map  $X'\to X$ where $X'$ is a smooth projective variety birational
to the product of $\PP^1$ and a smooth threefold. Therefore
$X'$ satisfies the Hodge conjecture by lemma \ref{lemma:fiber}
and corollary \ref{cor:dim5}. Consequently the conjecture holds
for $X$ by corollary \ref{cor:lifthodge2}.
\end{proof}

\begin{lemma}\label{lemma:torusaction} 
Suppose $X$ is a smooth projective variety on which an algebraic torus
 $T = (\C^*)^N$ acts. $X$ satisfies $GHC(H_i(X),p)$
if every component $S\subset X^{T}$ of the fixed point set
satisfies $GHC(H_i(S),p)$. In particular, $X$ satisfies the
 Hodge (respectively generalized Hodge) conjecture if the components of
the fixed point locus $X^{T}$ do.
\end{lemma}

\begin{proof}
Let $\C^*\to T$ be given by $t\mapsto (t^{a_1},t^{a_2},\ldots)$.
By choosing $1<< a_1 << a_2<<\ldots$, we can arrange equality of
the fixed point sets $X^{T} = X^{\C^*}$. Thus we can assume
that $T = \C^*$.
 By a theorem of  Bialynicki-Birula \cite{bb}, 
$X$ has a stratification such that the strata are
affine space  bundles over components of the fixed point locus.
The result follows from corollary \ref{cor:strat} and
lemma \ref{lemma:fiber}.
\end{proof}

\begin{cor} The (generalized) Hodge conjecture holds if
 $dim X^{T}\le 3$ ($\le 2$).
\end{cor}

When the fixed points are isolated,  Bialynicki-Birula's
decomposition is a cellular decomposition, so
all the cohomology is algebraic.

\section{Hodge conjecture for products}

In this section, we  give some criteria for the 
Hodge conjecture to hold for a power or self product of a variety  $X$.
The basic tool for this is the Mumford-Tate group of a Hodge structure.
The main reference is  \cite[I, sect. 3, sect. 5]{dm},
although the summary in \cite[II]{gordon} should be  sufficient and
a bit more accessible. Let  us start with an abstract characterization
since it explains the significance most clearly.
Given a collection of rational Hodge structures $H_i$, let $\langle H_i\rangle$ 
be the tensor $H_i$. Tannakian considerations
show that $\langle H, \Q(1)\rangle$ is equivalent to the category of rational
representations of a canonically determined 
algebraic group defined over $\Q$; this is the Mumford-Tate
group $MT(H)$. By definition,  sub-Hodge structures of $H^{\otimes n}$ are 
the same thing as $MT(H)$-submodules.
There is a homomorphism $\tau:MT(H)\to MT(0) =\mathbb{G}_m$
induced by the inclusion $\langle\Q(1)\rangle\subset \langle
H,\Q(1)\rangle$. 
The Hodge  
or special Mumford-Tate group $Hdg(H)$ is the kernel of this map.
The key property  which  makes this notion useful is that the
$Hdg(H)$ invariant tensors in $H^{\otimes n}$ are exactly
the Hodge cycles in this space ($MT(H)$ will act on the
space of Hodge cycles through a power of $\tau$).
An  alternative characterization is as follows: the real Hodge structure
$H\otimes \R$ has an action of $\C^*$ (viewed as a real group).
$Hdg(H)$  can be defined as  the
smallest $\Q$-algebraic subgroup of $GL(H)$ whose real points
contain the image of $U(1)\in GL(H\otimes \R)$.
 When $H$  has a polarization $\psi$ (which is a  symmetric or
 alternating form or according to the weight), then $Hdg(H)$ is a reductive 
subgroup of  $SO(\psi)$ or $Sp(\psi)$.

The best understood (nontrivial) case is that of an
abelian variety $X$. We  refer the reader to
the survey articles \cite{gordon, murty}
for further details and references.
In this case, let $Hdg(X)= Hdg(H^1(X))$.
Given a polarization $\psi$ of $X$,
the Lefschetz group $Lef(X)$,
is the centralizer of $End(X)\otimes \Q$ in $Sp(H^1(X),\psi)$.
The Lefschetz group turns out to be independent
of the polarization, and it always contains the Hodge group.

 Recall by Poincar\'e reducibility \cite{mumford}, any abelian
variety is isogenous to a product of simple abelian varieties.
Let $X$ be a simple abelian variety, then $E = End(X)\otimes \Q$ will
be a division algebra over its  center  $K$ which is a number field.
Albert [loc. cit., p 201] has classified the possiblities:

\begin{enumerate}
    \item[I] $E=K$ is a totally real field. 
    
    \item[II] $K$ is totally real, and $E$ is a totally indefinite
    quaternion algebra over $K$.

    \item[III] $K$ is totally real, and $E$ is a totally definite
    quaternion algebra over $K$.
    
    \item[IV] $K$ is a CM field.

 \end{enumerate}

\begin{thm}[Murty, Ribet]\label{thm:HdgLef}
If $X$ is an abelian variety, then
 $Hdg(X)=Lef(X)$ if and only if all Hodge classes on all
 powers of $X$ are products of divisor classes.
 In particular, equality of these groups implies
 that the Hodge conjecture holds for all powers
of $X$.
\end{thm}

\begin{remark}\label{remark:HdgLef}
There are a number of interesting cases where these conditions are
satisfied. $Hdg(X) = Lef(X)$ when:

\begin{enumerate}

\item
$X$ is a simple abelian variety with $dim \, X$ an odd 
prime (Tankeev, Ribet \cite{ribet}),

\item
$X$ an abelian variety of dimension 2 or 3 (this appears
to be part of the standard folklore; a proof
can be found in \cite{moonzar}).

\item $E=End(X)\otimes \Q$ is a totally real number field such that
$dim X/[E:\Q]$ is odd [loc. cit.)] (it follows that $X$ is simple of
type I),

\item $E$ is a CM field such that $dim X/[E:\Q]$ is prime
[loc. cit.] (in particular, $X$ is simple of type IV),

\item
$X$ is isogenous to a product of elliptic curves (see Murty \cite{murty2}),

\item
 $X$ is the Jacobian of a very general curve (corollary
 \ref{cor:vgjacobian}),

\item $X$ is a very general abelian variety. (This essentially 
goes back to Mattuck; a proof can be given along the
lines of that of corollary \ref{cor:vgjacobian}.)

\item The Jacobian of a modular curve 
$X_1(N)$ which is the compactification of the quotient of
the upper half plane by
$$\Gamma_1(N) = \{
\left(\begin{array}{cc}a& b\\ c& d\end{array}\right)
    \,|\, c\equiv 0, a \equiv 1\, (mod\> N)\}$$
(Hazama, Murty)

\end{enumerate}\end{remark}

\begin{cor}\label{cor:HdgLef}
 The Hodge conjecture holds for all powers of $X$, provided
 that it is on the above list.
\end{cor} 

Another class of examples where the conclusion of the corollary is
known, even though the Lefschetz and Hodge groups may differ,
 comes from the work of Shioda  \cite{shioda}. 

\begin{thm} If $J$ is the Jacobian of the Fermat curve $x^m + y^m +z^m=0$
where $m$ is prime or $m\le 20$, then the Hodge conjecture
holds for all powers of $J$,
\end{thm}

\begin{proof} This is an immediate consequence of \cite[4.4]{shioda}.
\end{proof}

Hazama \cite{hazama} has establised the generalized Hodge conjecture
in certain of
the above cases [loc. cit., pg 201]. The theorem can be formulated as follows:

\begin{thm}[Hazama] Let $X$ be an abelian variety satsifying $Hdg(X) = Lef(X)$
and such that all simple factors are of types I or II, then
the generalized Hodge conjecture holds for $X$.
\end{thm}

\begin{cor}\label{cor:Hazama}
If $X$ is as above, then the generalized Hodge conjecture holds for
all powers of $X$.
\end{cor}

\begin{proof} $Hdg(X^k) = Hdg(X)$
since $H^1(X)$ and $H^1(X^k) = H^1(X)^k$ generate the
same tensor category. Also $Lef(X) = Lef(X^k)$ \cite[cor. 4.7]{milne}.
Therefore $X^k$ satisfies the same conditions as the theorem.
\end{proof}

\begin{cor}\label{cor:ghc4curves}
 If  $E=End(X)\otimes \Q$ is a totally real number field such that
$dim X/[E:\Q]$ is odd
then the generalized Hodge conjecture holds for all powers $X$.
\end{cor}

\begin{proof} The conditions imply that $X$ is simple of type I.
Also by remark~\ref{remark:HdgLef}, $Hdg(X) = Lef(X)$.
\end{proof}

The corresponding results for curves follows from:

\begin{prop}\label{prop:jacobian}
The  Hodge (respectively generalized Hodge) conjecture holds for all powers of a smooth projective
curve $X$ if and only
if it holds for all powers of its Jacobian $J(X)$.
\end{prop}

\begin{proof}
Suppose that the (generalized) Hodge conjecture holds for $J(X)^k$.
Choose a  base point on $X$, and let $\alpha:X\to J(X)$ be the 
corresponding Abel-Jacobi map. Then $\alpha^*:H^*(J(X))\to H^*(X)$
is a surjection. The K\"unneth formula implies that the
induced map $\alpha^k:X^k\to J(X)^k$ also induces a surjection
on cohomology. Therefore the (generalized) Hodge conjecture holds for
$X^k$ by corollary \ref{cor:lifthodge}.

Let $g$ be the genus of $X$. Consider the map $\beta:X^{g}\to J(X)$ given
by $(x_1,\ldots x_g) \mapsto \alpha(x_1)+\ldots \alpha(x_g)$.
This map induces a surjection $\beta^k_*:H^*(X^{gk},\Q)\to H^*(J(X)^k,\Q)$
by lemma \ref{lemma:rathomology}. Thus
the (generalized) Hodge conjecture for $X^{gk}$
implies it for $J(X)^k$.
\end{proof}

Define $Hdg(X)= Hdg(J(X))$ and $Lef(X) = Lef(J(X))$. Then:

\begin{cor}\label{cor:Lef}
If $Hdg(X)=Lef(X)$ or if $X$ is Fermat of degree $m\le 20$ or a prime,
then the Hodge conjecture holds for all powers of $X$.
\end{cor}

The  characterization of Mumford-Tate groups \cite[p. 43]{dm}
together with  \cite[7.5]{deligne} (see also
\cite[2.2-2.3]{schoen}) yields:

\begin{lemma}\label{lemma:mono}
Given a polarized integral variation of Hodge structure $V$ over
a smooth irreducible complex variety $T$, there exists a
countable union of proper analytic subvarieties $S\subset T$
such that $Hdg(V_t)$
 contains a finite index subgroup of the monodromy group
$$image[\pi_1(S,t)\to GL(V_t)]$$
for $t\notin S$.
\end{lemma}

\begin{cor}\label{cor:vgjacobian}
If $X$ is
very general in the moduli space of curves, then  the generalized
 Hodge conjecture holds  for all powers of $X$.
\end{cor}

\begin{proof}
 Choose $n \ge 3$ and let $M_{g,n}$ be the fine moduli space
 of smooth projective curves of genus $g$ with level $n$ structure
\cite[13.4]{AO}.
Let  $\pi:{\mathcal X}\to M_{g,n}$ be the universal curve.  Lemma
 \ref{lemma:mono} applied to $\R^1\pi_*\Z$ shows that
 there exist a countable union of proper
subvarieties $S'\subset M_{g,n}(\C)$ such that a finite index subgroup of
the monodromy group 
$$\Gamma=image[\pi_1(M_{g,n}, t)\to  GL(H^1(\mathcal{X}_t))]$$
is contained in  $Hdg(\mathcal{X}_t)$ for each 
$t\notin S'$. Let $S$ be the image of $S'$ in $M_g(\C)$.
By Teichmuller theory, any finite index subgroup of $\Gamma$ is seen to be Zariski dense in 
the symplectic group (see \cite[12]{hain}). Hence the Hodge
group contains the symplectic group whenever $t\notin S$.
But this forces 
$$Hdg(\mathcal{X}_t)=Lef(\mathcal{X}_t) = Sp(H^1(\mathcal{X}_t)).$$

Fix $X= \mathcal{X}_t$, with $t$ as above.
We will show that 
$End(X)\otimes \Q = \Q$, and this  will  finish the proof
by corollary  \ref{cor:ghc4curves}.
The natural  map 
$$End(X)\otimes \Q \to  End(H^1(X,\Q))$$
is injective, and the image lies in the ring $End_{MHS}(H^1(X))$ 
of endomorphisms of the Hodge structure $H^1(X)$. This is
contained in  the space $Hdg(X)$-equivariant maps automorphisms
of  $H^1(X)$. $Hdg(X)$ acts irreducibly, since it is
the full symplectic group. Therefore Schur's lemma implies that
$End(X)\otimes \Q = \Q$ as claimed.
\end{proof}

Given a smooth projective surface $X$, $H^2(X,\Z)$ carries a
symmetric bilinear form $<,>$ given by cup product.
Let $A(X) = H^2_{alg}(X,\Q)$ or equivalently the Neron-Severi
group tensored with $\Q$. The rational transcendental lattice $T(X)$ 
 is the orthogonal complement 
$A(X)^{\perp}$. The decomposition
$$H^2(X,\Q) = A(X)\oplus T(X)$$
respects Hodge structure. 

\begin{lemma}\label{lemma:characterizationofT}
$T(X)$ is the smallest rational Hodge substructure containing $H^{20}(X)$.
\end{lemma}

\begin{proof}
  Let $V$ be the smallest Hodge  substructure containing
 $H^{20}(X)$. Then $V^{\perp}$ is a rational subspace lying
in $H^{11}(X)$. Therefore $V^{\perp}\subseteq A(X)$ by the 
Lefschetz $(1,1)$ theorem, and this implies
$T\subseteq V$. On the other hand, the Hodge-Riemann bilinear
relations imply that $H^{20}(X)$ is orthogonal to $H^{11}(X)$.
Therefore $H^{20}(X)\subseteq T(X)$, and this give the opposite
inclusion.
\end{proof}

If $Z\subset Y\times X$ is a codimension $2$ correspondence, then the
induced morphism $Z_*:H^i(Y,\Q)\to H^i(X,\Q)$ can be restricted
to $H^{i0}$ to get a map
$H^0(\Omega_Y^i)\to H^0(\Omega_X^i)$. This  can
be interpreted directly in terms of differential forms.
After, replacing $Z$ by a resolution $\tilde Z$, the map is a composition of
the pullback $H^0(\Omega_Y^i)\to H^0(\Omega_{\tilde Z}^i)$ and the trace 
$H^0(\Omega_{\tilde Z}^i)\to H^0(\Omega_X^i)$. By definition $Z_*$ preserves
$A$, and it preserves $T$ by lemma~\ref{lemma:characterizationofT}.

\begin{prop}\label{prop:transclat}
Let $X$ and $Y$ be smooth projective surfaces, and
let $Z\subset Y\times X$ be a codimension $2$ correspondence.
 The maps
 $$Z_*:H^0(Y,\Omega_Y^i)\to H^0(X,\Omega_X^i),\> i=1,2$$
are surjective if and only if
$$Z_*:H^1(Y)\to H^1(X)
\mbox{ and }
Z_*:T(Y)\to T(X)$$
are surjective.
Suppose that  these  conditions hold and that
 $Y^k$ satisfies the  Hodge (respectively generalized Hodge) conjecture for
all $k\le n$,
then $X^k$ satisfies the  Hodge (respectively generalized Hodge) conjecture for all $k\le n$. 
\end{prop}

Before giving the proof, note that the symmetric group $S_n$ acts
on $X^n$, and hence on $H^*(X^n)$  by permutation of factors.
This action is compatible with the K\"unneth decomposition in
the following sense: given classes $\alpha_{j}\in H^{*}(X)$
$$\sigma(\alpha_{1}\times \alpha_{2}\times \ldots \alpha_{n})= 
\alpha_{\sigma(i_1)}\times \alpha_{\sigma(i_2)}\times \ldots \alpha_{\sigma(i_n)}
$$

\begin{proof}
The equivalence of the surjectivity statements follows
from elementary Hodge theory and lemma~\ref{lemma:characterizationofT}. 
    It is enough to prove the remainder of the proposition for $k=n$.
 The correspondence $Z^n\subset Y^n\times X^n$
induces morphisms $H^i(Y^n)\to H^i(X^n)$. 
This is compatible with the K\"unneth decompositions
$$H^{i}(Y^n) = \bigoplus_{j_1,j_2,i_1\ldots,\sigma}
\sigma[N(Y)^{\otimes j_1}\otimes T(Y)^{\otimes j_2}\otimes H^{i_1}(Y)\otimes \ldots
H^{i_{n-1}}(Y)]$$
$$H^{i}(X^n) = \bigoplus_{j_1,j_2,i_1\ldots,\sigma}
\sigma[N(X)^{\otimes j_1}\otimes T(X)^{\otimes j_2}\otimes H^{i_1}(X)\otimes \ldots
H^{i_{n-1}}(X)],$$
where  $i_\ell\not= 2$,   $2j_1+2j_2+\sum i_\ell=i$  
and $\sigma$ ranges over a set of permutations.
By hypothesis, any Hodge cycle in 
$$N(X)^{\otimes j_1}\otimes T(X)^{\otimes j_2}\otimes
H^{i_1}(X)\otimes \ldots$$
can be lifted to a product of an algebraic cycle with a
Hodge cycle on $H^{i-2j_1}(Y^{n-j_1})$...
\end{proof}

\begin{prop}\label{prop:hodgesurface}
Let $X$ be a smooth projective surface, then:

\begin{enumerate}
 
\item[a)] If $f:Y\dashrightarrow X$ is a dominant rational map of smooth
projective surfaces,
$X^k$ satisfies the Hodge (respectively generalized Hodge) conjecture  for all $k\le n$ 
if $Y^k$ does for all $k\le n$.
In particular, the condition is birationally invariant.

\item[b)] If $X$ is a rational surface, then $X^k$ satisfies the
generalized Hodge conjecture for all $k>0$.

\item[c)] If $X\to C$ is a (possibly nonminimal) ruled surface, 
then $X^k$ satisfies the  Hodge (respectively generalized Hodge)  conjecture for all $k\le n$ if 
it holds for all $C^k$ in the same range. In particular, the
Hodge conjecture holds for all $k$ if $Hdg(C) =Lef(C)$ or $C$ is 
Fermat of degree$\le 20$ or a prime.

\item[d)] If $X$ is an abelian surface, then $X^k$
 satisfies the Hodge conjecture for all $k$.

\item[e)] If $X$ is a K3 surface such that $(T(X), <,>)$ can be embedded
  isometrically in $H^{\oplus 3}$, where $H = \Q^2$ with the quadratic
form $\left(\begin{array}{cc}0 & 1\\ 1& 0\end{array}\right)$, then  $X^k$
 satisfies the Hodge conjecture for all $k$. 

\item[f)] If $X$ is an Kummer surface or a K3 surface with Picard
  number $\ge 19$, then $X^k$
 satisfies the Hodge conjecture for all $k$.
\end{enumerate}
\end{prop}

\begin{proof}

a) Let $f:Z\to Y$ be a sequence of blow ups such that $Y\dashrightarrow X$
extends to a morphism $g:Z\to X$. $f^*$ induces  isomorphisms
$H^1(Y)\to H^1(Z)$ and $T(Y)\cong T(Z)$, and $g_*$ induces surjections
$H^1(Z)\to H^1(X)$ and $T(Z)\to T(X)$.
Therefore a) follows from proposition~\ref{prop:transclat}.




b) By a), we can assume that $X=\PP^2$. The statement is clear
since $X^n$ has a cellular decomposition.

c) By a), we can assume that $X$ is minimal. The statement
follows from lemma \ref{lemma:fiber} and corollary \ref{cor:Lef}.

d) $X$ is either simple or isogenous to a product of
two elliptic curves. So the result follows from
corollary \ref{cor:HdgLef}.

e) By a theorem of Mukai \cite[1.12]{mukai2}, there is an
abelian surface $A$ and a correspondence on $X\times A$ inducing
an isomorphism $T(A)\cong T(X)$. Therefore the result follows
from d) and proposition \ref{prop:transclat}.

f) Follows from e) since the conditions imply that $T(X)$ embeds
into $H^3$.

\end{proof}

Finally, we give some criteria for the Hodge conjecture to
hold when the Hodge group is large.

\begin{thm}\label{thm:oddpowers} Suppose that $Y$ is a smooth
   projective variety
  such that $d= dim\, Y$ is odd and
\begin{enumerate}
\item $H^i(Y) = 0$ when  $i$ is odd and different from $d$
\item $H^i(Y)=H^i_{alg}(Y)$ when $i$ is even.
\item  $Hdg(H^d(Y))$ coincides with the symplectic
group $Sp(H^d(Y))$ with respect to the cup product.

\end{enumerate}
Then the Hodge conjecture holds for all powers of $Y$.
\end{thm}

A somewhat weaker analogue can be proven for even dimensional varieties.
As noted earlier, the symmetric group $S_m$ acts on $H^i(Y^m)$.
Call an element $\alpha\in H^i(Y^m)$ antisymmetric if $\sigma(\alpha)
= (-1)^\sigma \alpha$, where $(-1)^\sigma$ denotes the parity
of $\sigma$.

\begin{thm}\label{thm:evenpowers} Suppose that $Y$ is a smooth
  projective variety
  such that $d= dim\, Y$ is even and
\begin{enumerate}
\item $H^i(Y) = 0$ when  $i$ is odd
\item $H^i(Y)=H^i_{alg}(Y)$ when $i$ is even and different from $d$.
\item   $SO(T)\subseteq Hdg(H^d(Y))$ for 
 $A \subseteq H^{d}_{alg}(Y)$ (equality is not required) 
and $T = A^{\perp}\subset H^d(Y)$.

\end{enumerate}
Then there exists an antisymmetric Hodge class $\delta \in H^{d\tau}(Y^\tau)$,
where $\tau = dim T$, such that any Hodge class on a $Y^n$
is a linear combination $\alpha+\sum\sigma_i(\delta\times \beta_i)$
 where $\alpha$ and $\beta_i$ are algebraic cycles, and $\sigma_i\in S_n$.
\end{thm}

The proof of these theorems is an exercise in invariant
theory. Let $V$ be a $d$-dimensional vector space over a field 
of characteristic $0$ with a nondegenerate alternating  or
symmetric bilinear form $\psi$.
Let $G = Sp(V,\psi)$ in the first case, and 
let $G=SO(V,\psi)$ in the second.
The form induces an isomorphism $V\cong V^*$ as $G$ modules.
Therefore
$\psi$ can be regarded as a tensor in $V\otimes V$ which is invariant
under the $G$. Choose a nonzero element $\delta\in \wedge^dV$
which we identify with an antisymmetric tensor in $V^{\otimes d}$.
This is also $G$-invariant.  The symmetric
group $S_{N}$ acts on $V^{\otimes N}$ by permuting factors, and this
action commutes with the $G$ action. Tensor
products of the previous tensors, and their transforms under the symmetric
group, generate all  $G$-invariant tensors. More precisely:

\begin{lemma}[Weyl] \label{lemma:inv}
With the above notation
\begin{enumerate}
\item $V^{\otimes n}$ is spanned by the $S_n$ orbit of
  $\{\psi^{\otimes n/2}\}$ if $G = Sp(V,\psi)$, and
\item $V^{\otimes n}$ is spanned by the union of the $S_n$ orbits of
  $\{\psi^{\otimes n/2}\}$ and $\{\delta\otimes \psi^{\otimes
  (n-d)/2}\}$ otherwise;
\end{enumerate}
where these sets are taken to be empty unless  the exponents are
nonnegative integers.
\end{lemma}

\begin{proof} See \cite[F.13, F.15]{FH}
\end{proof}

\begin{proof}[Proof of Theorems \ref{thm:oddpowers} and \ref{thm:evenpowers}]
The proofs of both theorems will run mostly in parallel with
occasional branching. In the case of theorem \ref{thm:evenpowers}
 we can assume that $dim\, T>1$, otherwise
the theorem is vacuous. Since $T$ can be identified with
$H^d(Y)/A$, it  carries a natural Hodge  structure.
 Set $V= H^d(Y)$.
If $d$ is odd, let $\psi$ denote the cup product
form on $V$, and let $G = Sp(V,\psi)$.
If $d$ is even, let $\psi$ denote the restriction of
the cup product form to $T$, and let $G=SO(T,\psi)$.
By hypothesis $Hdg(V)$ contains $G$.
 Therefore,
the Hodge cycles in $V^{\otimes i}(j)$ 
are $G$-invariant. Since $T$ is irreducible and nontrivial,
this proves theorem \ref{thm:evenpowers} for the first power $Y^1$.

We will show that the class $\psi$ is represented
by an algebraic cycle on $Y\times Y$.
In this paragraph we treat the case where $d$ is odd.
By K\"unneth's formula we have
$$H^{2d}(Y\times Y) = [V\otimes V] \oplus W$$
where
$$W=\bigoplus_{i\not= d}[H^i(Y)\otimes H^{2d-i}(Y)].$$
Then by the first two assumptions of theorem \ref{thm:oddpowers}, 
$W\subset H^{2d}_{alg}(Y\times Y)$. In particular, the sum $\Delta'$ of the
K\"unneth components of the diagonal $\Delta$ in $W$ is algebraic.
By lemma \ref{lemma:inv}, $[V\otimes V]^{G}$ is one dimensional.
Therefore it must be spanned by $\Delta-\Delta'$.
In particular, the  form $\psi$ in $V\otimes V$ is  algebraic.

Now suppose $d$ is even and $dim\, T\ge2$. We decompose $V\otimes V$ further as
$[T\otimes T] \oplus V_1\oplus V_2$ where
$V_1 = [T\otimes A] \oplus [A\otimes T]$ and $V_2 = A\otimes A$.
$V_1$ is isomorphic to a sum of copies of Tate twists of $T$.
Therefore, for example by lemma \ref{lemma:inv}, $V_1^G=0$.
This implies that the components of  $\Delta$  
in $V_1$ are $0$.
We have $W\oplus V_2\subset H^{2d}_{alg}(Y\times Y)$, therefore
the sum $\Delta'$ of the components of
 $\Delta$ in $W\oplus V_2$ is algebraic.
By lemma \ref{lemma:inv}, $[T\otimes T]^G$ is one dimensional.
Therefore $\Delta-\Delta'$ spans it, and this
proves algebraicity of $\psi$ in this case.

Let $m >0$.
The K\"unneth decomposition can be written as
$$H^{i}(Y^m) = \bigoplus_{i_k,\sigma}
\sigma[H^d(Y)^{\otimes j}\otimes H^{i_1}(Y)\otimes \ldots
H^{i_{m-1}}(Y)],$$
where each $i_k\not= d$, $dj=i-\sum i_k$  and $\sigma$ ranges over a set of permutations.
When $d$ is even, we refine this further as
$$H^{i}(Y^m) = \bigoplus_{i_k,j_\ell,\sigma}
\sigma[T^{\otimes j_1}\otimes A^{\otimes j_2}\otimes H^{i_1}(Y)\otimes \ldots
H^{i_{m-1}}(Y)],$$
for an appropriate set of indices.
Let $\delta$ be the nonzero rational element of 
$\wedge^\tau T\subset H^{d\tau}(Y^\tau)$; it is necessarily an
antisymmetric Hodge cycle.
By our assumptions, we can choose bases $\{\gamma_l^{(i_k)}\}$ for the groups $H^{i_k}(Y)$
(and $\{\gamma_l'\}$ for $A$ if $d$ is even)
consisting of algebraic cycles. Then, we have
$$H^{i}(Y^m) = \bigoplus_{\Gamma,\sigma} \sigma[H^d(Y)^{\otimes j}\times \Gamma] $$
if $d$ is odd, or
$$H^{i}(Y^m) = \bigoplus_{\Gamma,\sigma} \sigma[T^{\otimes j}\times \Gamma] $$
if $d$ is even,
where $\Gamma$ ranges over products of the $\gamma$'s (and $\gamma'$'s).
A Hodge cycle on $Y^m$ can be decomposed into a sum of Hodge
cycles from each of the above summands.
As an abstract Hodge structure  each of these
summands is isomorphic to a Tate twist of $ H^d(Y)^{\otimes j}$ or 
$T^{\otimes j}$.
As already noted the Hodge  cycles here are $G$-invariant.
By lemma \ref{lemma:inv}, the  $G$-invariant classes
are spanned by the orbit of $\psi\times\psi\times\ldots \psi$
and when $d$ is odd.
These elements
are algebraic by the previous paragraph, and this finishes the proof
of theorem~\ref{thm:oddpowers}. 
When $d$ is even, we have
additional $G$-invariant classes lying in the orbit
of  $\delta\times \psi\ldots\psi$, and this accounts for the
remaining terms in the statement of theorem~\ref{thm:evenpowers}.
\end{proof}

\begin{cor}\label{cor:symevenpowers}
Suppose that $Y$ satisfies the assumptions of 
theorem~\ref{thm:evenpowers}. Then 
the Hodge conjecture holds for arbitrary symmetric powers $S^nY$.
\end{cor}

\begin{proof}
We may assume that $T\not= 0$, or else there is nothing to
prove. In this case, $dim T \ge 2$, since $T$ would contain
nontrivial elements of type $(p,q)$ and $(q,p)$ with $p\not= q$.

Consider the map $Y^n\to S^nY$. It induces a surjection
$H^i(Y^n)\to H^i(S^nY)$ by lemma~\ref{lemma:rathomology}.
Therefore a Hodge cycle $\gamma$ on $S^nY$ can be lifted
to a Hodge cycle $\gamma'$ on $Y^n$. $\gamma'$ can be expressed
as a sum  $\alpha+\sigma_1(\delta\times \beta_1)+\ldots$ as in 
theorem~\ref{thm:evenpowers}.
Since $S_n$ acts trivially on $S^nY$, the map $H^i(Y^n)\to H^i(S^nY)$
factors through the space of coinvariants $H^i(Y^n)_{S_n}$.
The image of  $\sigma_k(\delta\times \beta_k)$ in $H^i(Y^n)_{S_n}$
must vanish because of the antisymmetry of $\delta$. (More explicitly:
 $\sigma_k(\delta\times \beta_k)$, $\delta\times \beta_k$
and $(12)(\delta\times \beta_k)= -\delta\times \beta_k$ have the
same image.)
Therefore $\gamma$ is represented by the image of $\alpha$.
\end{proof}

\begin{cor}\label{cor:sympowersofK3}
Let $Y$ be a smooth projective surface with $q=0$ and  $p_g=1$
(e. g. a K3 surface)
such that $End_{MHS}(T(Y))=\Q$, then the Hodge conjecture
holds for all symmetric powers of $Y$.
\end{cor}

\begin{proof}
By \cite[2.2.1]{zarhin2} the hypothesis implies that
$Hdg(T(X))=SO(T(X))$.
\end{proof}

\begin{cor}\label{cor:hyp1}
Suppose that $Z\subseteq \PP^N$ is an even dimensional  smooth
projective variety with a cellular decomposition.  If $H$
is a sufficiently general hyperplane section of $Z$, then the Hodge 
conjecture holds for all powers of $H$.
\end{cor}

The precise meaning of ``sufficiently general'' above is the following:
For any Lefschetz pencil  $\mathcal{Y}\to \PP^1$  of hyperplane
sections of $Z$, $H$ can be taken to be $\mathcal{Y}_t$ for all but countably
many $t\in \PP^1(\C)$.

\begin{proof}
Let $dim Z = d+1$.
Conditions 1 and 2 of  theorem \ref{thm:oddpowers} hold for any smooth
hyperplane section  by the
weak Lefschetz theorem ($i<d$) and the hard Lefschetz theorem ($i>d$).
Suppose that  $\mathcal{Y}\to \PP^1$ is a Lefschetz pencil,
and let $U\subset \PP^1$ be the complement of the set of
critical values.
Then for any smooth fiber $Y_t$, 
$H^d(Y_t)=H^d(Z)\oplus E = E$, where $E$ is the subspace generated
by vanishing cycles \cite[5.5]{sga7}. Since the orthogonal complement
$E^\perp =0$, we can apply 
the Kazhdan-Margulis theorem \cite{deligne2} to see that the image
of the monodromy representation $\pi_1(U,t)\to GL(H^d(Y_t))$
is Zariski dense in $Sp(H^d(Y_t))$.
Then by lemma \ref{lemma:mono}
there exists a countable set $S$
such that if $t\notin S$ then $Y_t$ satisfies the third condition
of theorem \ref{thm:oddpowers}.
\end{proof}

It is possible to prove a  weaker statement
when $dim \, Z$ is odd.

\begin{cor}\label{cor:hyp2}
Suppose that $Z\subseteq \PP^N$ is a smooth projective variety
with a cellular decomposition such that $d=dim Z-1$ is even. 
There exists  an integer $n_0$  such that if $H$
is a sufficiently general hypersurface section of $Z$ of degree $n\ge n_0$, 
then the Hodge conjecture holds for all symmetric powers $S^mH$.
\end{cor}

\begin{proof}
The proof of this is
very similar to that of corollary \ref{cor:hyp1}.
Suppose that  $\mathcal{Y}\to \PP^1$ is a Lefschetz pencil
of hypersurfaces of degree $n$,
and let $U\subset \PP^1$ be the complement of the set of
critical values.
Then for any smooth fiber $Y_t$, there is an orthogonal decomposition
$H^d(Y_t) = A\oplus T$ where $A = H^d(Z)$ and $T$ is the  subspace generated by
 vanishing cycles \cite[5.5]{sga7}.
 The monodromy representation $\pi_1(U,t)\to
GL(T)$ is dense in $O(T)$ for $n$ greater than
or equal to some $n_0>0$ by \cite[thm B]{verdier}.
Therefore any finite index subgroup $\Gamma\subseteq \pi_1(U,t)$
contains a subgroup of finite index which dense in  $SO(T)$.
Then by lemma \ref{lemma:mono}
there exists a countable set $S$
such that if $t\notin S$ then $Y_t$ satisfies the third condition of
theorem \ref{thm:evenpowers},
and the first two conditions are automatic. Therefore the result
follows from corollary~\ref{cor:symevenpowers}.
\end{proof}

\section{Moduli of vector bundles  over curves}

When  $X$ is a smooth projective curve,
 let $U_X(n,d)$ be the
 moduli space of semistable bundles on $X$ of rank $n$ and 
degree $d$. It is smooth and projective if $n$ and 
$d$ are coprime.

\begin{thm}\label{thm:powersXtoM} Let
$X$ be a smooth projective  curve and  $M=U_X(n,d)$ with $n$ and 
$d$ coprime. If the  Hodge (respectively generalized Hodge) conjecture holds for
all powers $X^k$, then the Hodge (respectively
generalized Hodge) conjecture holds for $M$.
\end{thm}

\begin{cor}[del Ba\~no \cite{delb}]
If the Hodge conjecture holds for all powers of $J(X)$, then it
holds for $M$.
\end{cor}

\begin{proof}
Proposition \ref{prop:jacobian}.
\end{proof}

We record specific instances where this holds:

\begin{cor} If $X$ is
\begin{enumerate}

\item
very general in the moduli space of curves (Biswas-Narasimhan \cite{BN}),  
\item
a curve of genus 2 or 3,
\item
a curve of prime genus such that the Jacobian is simple, or
\item
a Fermat curve $x^m + y^m +z^m = 0$ with $m$ prime or less than $21$, 
or,
\item a curve admitting a surjection from an  $X_1(N)$,
\end{enumerate}
then the Hodge conjecture holds 
for $M$ (as above).
\end{cor}

\begin{proof}
These follow from the results of section 3. For the second case,
the Jacobian is either simple or isogenous to a product  of
elliptic curves. For the last case, $X^k$  satisfies the Hodge
conjecture since it is dominated by $X_1(N)^k$.
\end{proof}

The conclusion of the first part can be strengthened considerably:

\begin{cor}
If $X$ is very general in the moduli space of curves, then 
the generalized Hodge conjecture holds for $M$.
\end{cor}

\begin{proof}
This follows from corollary \ref{cor:ghc4curves} and 
proposition \ref{prop:jacobian}.
\end{proof}

We give two proofs of the theorem. The first is fairly elementary,
while the second is easier to generalize.
Let $E$ be a vector bundle on $X$ and $Q_n(E)$ denote the
the ``Quot'' scheme parameterizing coherent subsheaves of $F\subset E$
such that $E/F$ is a torsion sheaf of length $n$. By \cite{groth}
 $Q_n(E)$ exists and is projective, and it can be seen to
be smooth since  $X$ is a curve.
When $E=O_X$, this is just the space of anti-effective divisors of 
degree $-n$. 
 
\begin{proof}[First proof]
Choose a divisor $D$ with $deg\, D >> 0$ (the precise requirements
will be given below). Let $E= O(D)^{\oplus n}$, $m= ndeg(D) -d$,
and $Q = Q_m(E)$. A point of $Q$ is given by an abstract vector
bundle $F$ on $X$ of rank $n$ and degree $d$ together with an
embedding $F\subset E$. Let $Q^s$ denote
the open subset which parameterizes pairs $(F\subset E)$ with $F$ stable.
There is an obvious map $\pi:Q^s\to M$.
 Choose $deg\, D$ sufficiently large,
so that $F^*(D)$ is globally generated and $H^1(F^*(D)) = 0$
for any stable vector bundle
$F\in M$. Thus any $F\in M$ can be embedded into $E$ which implies that
 $\pi:Q^s\to M$ is surjective. 
Let $\F$ denote a Poincar\'e bundle on $X\times M$. By our
assumptions,
$V = p_{X*}\mathcal{H}om(\F,p_M^*E)$ is locally free and commutes
with base change, where $p_X,p_M$ denote the projections of $X\times M$
to its factors.
Let $\pi':\PP (V^*)\to M$ denote the bundle of lines in the
fibers of $V$.
We have an embedding $Q^s\hookrightarrow \PP (V^*)$, which sends the
point $(F\subset E)$ to the line generated by the  corresponding
element of $Hom(F, E)$. The codimension of the complement
$Z= \PP(V^*)-Q^s$ is greater than or equal to $deg D$ \cite[8.2]{bgl}.
Let us assume that $deg D>(dim M+1)/2$.

The torus $T = (\C^*)^n$ acts on $E= O(D)^{\oplus n}$
by homotheties on each factor, and this induces an action on
$Q$ \cite{bgl}. The components of $Q^T$ are products of symmetric
powers of $X$ [loc. cit.]. In particular, the (generalized) 
Hodge conjecture holds for these components  by 
corollary~\ref{cor:lifthodge2}. Therefore the (generalized) Hodge
conjecture holds for $Q$ by lemma \ref{lemma:torusaction},
and consequently for $Q^s$ by proposition~\ref{prop:jannsen}.
We have an exact sequence
$$\Hbm_{2q-i}(Z)\to \Hbm_{2q-i}(\PP(V^*))\to\Hbm_{2q-i}(Q^s)\to\ldots$$
where $q = dim Q$. By our assumptions, 
$$\Hbm_{2q-i}(\PP(V^*))\cong \Hbm_{2q-i}(Q^s)$$
if $i \le 2dim M$. Therefore $GHC(H_{2q-i}(\PP(V^*)), 2(2q-i))$
(or $GHC(H_{2q-i}(\PP(V^*)), *)$) holds, or equivalently
 $GHC(H^i(\PP(V^*)), 2i)$
(or $GHC(H^i(\PP(V^*)), *)$) holds for $i \le 2dim M$.
The theorem now follows from corollary~\ref{cor:lifthodge}.
\end{proof}    

Embedded in this argument is a proof that $Q_n(E)$ satisfies the Hodge
conjecture for  $E=O(D)^n$. This is true more  generally:

\begin{prop}\label{prop:hc4Qn}
 If the  Hodge  (respectively generalized Hodge) conjecture holds for
all powers $X^k$, then the  Hodge (respectively generalized Hodge)  conjecture 
holds for $Q_n(E)$ for any $n$ and vector bundle $E$.
\end{prop}

The proof will be deferred to the next section.

Suppose $Y$ and $N$ are compact oriented manifolds and $c\in H^*(Y\times N,\Q)$.
We can decompose $c$ as
$$c=\sum p_Y^*d_i\cup p_N^*e_i=\sum d_i\times e_i$$
by the K\"unneth formula. The $e_i$ will be called  the 
K\"unneth components of $c$ along $N$. By Poincar\'e duality $c$ can 
be identified with a homomorphism $H^*(Y)\to H^*(M)$. Explicitly,
this is
$$c(d') = \sum \int_Y d_i\cup d'\, e_i$$
Thus the image of $c$ is contained in, and in fact equal to,
the span of the K\"unneth components.
The second proof will be based on
the following:

\begin{prop}\label{prop:ABtype}
    Suppose that $Y$ and $N$ are smooth projective varieties
    such that $Y^k$ satisfies the  Hodge (respectively generalized Hodge)  conjecture
    for all $k\le dim N$. In addition, assume that there exists
    a finite collection of algebraic correspondences 
    on $Y\times N$ such their K\"unneth components generate
    the cohomology ring $H^*(N,\Q)$. Then the  Hodge
    (respectively generalized Hodge) conjecture holds for $N$.
 \end{prop}   
 
 \begin{proof}
 Let $c_{j,i}\in H^{2i}(Y\times N)$ denote the cohomology classes
 of the above algebraic correspondences.
Taking exterior products, we get correspondences
$c_{j_1,i_1}(E)\times \ldots c_{j_n,i_n}(E)$ on $Y^n\times N^n$.
These can be 
pulled back along the diagonal map $N\hookrightarrow N^n$ to get
correspondences 
$C_{J,I}=C_{j_1,\ldots j_n,i_1\ldots i_n}$ on $Y^n\times N$. 
which induce morphisms (as in corollary \ref{cor:corr})
$$H^*(Y^n)(-(i_1+\ldots i_n-n))\to H^{*+2(i_1+\ldots i_n-n)}(N).$$
Let
$$A: \oplus_{J,I} H^{*}(Y^*)(*) \to H^a(N)$$
be the  sum of the $C_{J,I}$ maps as $I=(i_1,\ldots i_n)$
varies over all finite sequences with $\sum i_k\le a/2$ and $i_k > 0$.
The hypothesis implies that $H^a(N)$ is spanned by monomials in the
K\"unneth components of the $c_{j,i}$, and this is equivalent
to the surjectivity of the  map $A$.
Therefore corollary \ref{cor:lifthodge}
finishes the proof.
\end{proof}

 \begin{proof}[Second proof of theorem~\ref{thm:powersXtoM}]
Let $E$ be a Poincar\'e bundle on $X\times M$. The Chern classes give
correspondences $c_i(E)\in H^{2i}_{alg}(X\times M)$. 
The  K\"unneth components of these classes generate the cohomology
ring $H^*(M)$ by a theorem of Atiyah and Bott \cite[9.11]{AB} or
 more specifically Beauville's version of this theorem
 \cite{beauville}.
Thus we can apply the previous proposition.
\end{proof}
   
\begin{remark}
This proof actually gives a slightly stronger result that the
(generalized) Hodge conjecture holds for $M$ if it it holds for all powers
$X^k$ with $k\le dim M$.
\end{remark}

\section{Moduli of sheaves over surfaces}

For surfaces the analogous results are a bit more elusive.
We begin with the moduli space of ideals of zero dimensional
subschemes, i.e. the Hilbert scheme of points.  
We review some basic facts about these spaces; further
details can be found in \cite{gott}.
Let $X$ be a smooth projective variety. For each integer $n > 0$, 
 let $X^{(n)} = S^nX$ be the $n$th
symmetric power, and let $X^{[n]}$ be the Hilbert scheme of zero
dimensional subschemes of $X$ of length $n$.  There are canonical morphisms
$p:X^n\to X^{(n)}$ and $\psi:X^{[n]}_{red}\to X^{(n)}$. The
map $\psi$, called the Hilbert-Chow morphism, is birational.

\begin{thm}[Forgarty]
If $dim\, X= 2$, then $X^{[n]}$ is smooth (and projective)
for each $n$. It follows that $\psi:X^{[n]}\to X^{(n)}$
is a resolution of singularities.
\end{thm}

These spaces have a natural stratification.
Given a partition $\lambda = (n_1,n_2, \ldots n_k)$
of $n$ (i.e. a nonstrictly  decreasing sequence of positive
integers summing to $n$), let 
$$X^{(n)}_\lambda 
=\{p(x_1,\ldots x_n)\,|\,
 x_1= x_2=\ldots = x_{n_1}\not= x_{n_1+1} = \ldots = x_{n_1+n_2}\not=
 \ldots\}
$$
and let
$$X^{[n]}_\lambda = \psi^{-1} X^{(n)}_\lambda . $$
These are locally closed subsets  of $X^{(n)}$
and $X^{[n]}$ which will be regarded as subschemes with reduced
structure.
The scheme $X^{[n]}_{(n)}$ parameterizes $0$-dimensional subschemes
with support at a single point. There is a  morphism $\pi_n:X^{[n]}_{(n)}\to X$
which sends a subscheme to its support.

\begin{lemma}\label{lemma:gott1}(\cite[2.1.4, 2.2.4]{gott}) 
 $\pi_n$ is a locally
trivial fiber bundle. When $dim X = 2$, the fiber is smooth, projective
and has an  cellular decomposition.
\end{lemma}

Let $U_k\subset X^k$ be the open subset of $k$-tuples with
distinct components. For a partition $\lambda= (n_1,\ldots n_k)$
of $n$, define 
$$X_\lambda^{<n>} =U_{k}\times_{X^k}\prod_{i=1}^k X^{[n_i]}_{(n_i)}$$

\begin{lemma}\label{lemma:gott2}(\cite[2.3.3]{gott})
 $\ X^{[n]}_\lambda$ is a quotient of $X_\lambda^{<n>}$
by a finite group.
\end{lemma}

\begin{thm}\label{thm:hilb}
Let $X$ be a smooth projective surface such that all powers $X^k$,
with $k\le n$, satisfy the  Hodge  (respectively generalized Hodge) conjecture.
Then $X^{[n]}$ satisfies the  Hodge (respectively generalized Hodge) conjecture.
\end{thm}

\begin{proof}
Let  $\lambda = (n_1,\ldots n_k)$ be a partition of $n$.
By lemma \ref{lemma:gott1}, the $\pi_{n_i}$ are locally trivial
fiber bundles such that the fibers are smooth and projective
with cellular decompositions. Therefore 
$\prod  X^{[n_i]}_{(n_i)}$ is a fiber bundle over $X^k$
with these properties. Then lemma \ref{lemma:fiber} and  remark
\ref{remark:cell} imply that $\prod  X^{[n_i]}_{(n_i)}$
satisfies the (generalized) Hodge conjecture. Since  $X^{<n>}_\lambda$
is an open subset of this space, it also satisfies the
 (generalized) Hodge conjecture by proposition \ref{prop:jannsen}.
 
Lemma \ref{lemma:gott2} implies that $ X^{[n]}_\lambda$  has 
finite quotient singularities and that $X^{<n>}_\lambda$
is maps onto it.
Corollary \ref{cor:lifthodge2}, and the remarks preceding
it, yields the  (generalized) Hodge conjecture for
$\ X^{[n]}_\lambda$. Corollary \ref{cor:strat} finishes the
proof.
\end{proof}

\begin{cor} When $X$ and $k$ are as in proposition 
\ref{prop:hodgesurface}, $X^{[k]}$ satisfies the
Hodge conjecture.
\end{cor}



Proposition~\ref{prop:hc4Qn} can be proved by modifying the
above argument. We will sketch this.

\begin{proof}[Proof of proposition~\ref{prop:hc4Qn}]
 There is an analogue of the Hilbert-Chow morphism
$\psi:Q_n(E)\to S^nX=X^{(n)}$ which sends 
$$(F\subset E)\mapsto \sum_{p\in X} length(E_p/F_p)\, p$$
$X^{(n)}$ can be stratified as above. The strata satisfy
the Hodge conjecture since their closures are dominated
by powers of $X$.
The restriction
of $\psi$ to these strata are locally trivial fiber bundles
where the fibers are products of Grassmanians. Thus
the preimages of the strata under $\psi$ satisfy the Hodge
conjecture. We can now conclude the proof by
corollary \ref{cor:strat}.
\end{proof}

Let us turn to the higher rank case.
Given a smooth projective surface $X$, choose a 
polarization $H$ and elements $r\in \N$,
$c_i\in H^{2i}(X,\Z)$. Let
$$ch = r + c_1 + \frac{1}{2}(c_1^2-2c_2)$$
be the ``Chern character'', and  let $K$, $[X]$ and $td(X)$
respectively denote
the canonical, fundamental and Todd classes of $X$.
Call a class $c\in H^*(X,\Z)$ primitive if it is
 not multiple of a class  other than $\pm c$.
Let $M_X(r,c_1, c_2, H)$ be 
the moduli space of  torsion free sheaves
on $X$ of rank $r$ with Chern classes given by $c_i$  
which are $H$-semistable in the sense of
Gieseker-Maruyama \cite{gieseker} \cite{maruyama}. 
This is a projective scheme of finte
type. The open set of stable sheaves $M_X^s(r,c_1, c_2, H)$
tends to have more manageable local properties.
In particular, it is known to be smooth if $-K_X.H<0$
\cite[6.7.3]{maruyama}, or if $X$ is abelian or K3
\cite{mukai}.
In ideal cases, the stable locus coincides with the whole moduli space. 
Part of the standard folklore is:

\begin{lemma}\label{lemma:gcdofrc1chi}
If  $gcd(r, c_1\cdot H, ch\cdot td(X))=1$,  then
$M_X^s(r,c_1, c_2, H)=M_X(r,c_1, c_2, H)$
\end{lemma}

\begin{proof}
If $E$ lies in $M-M^s$, then
there exists subsheaf $F\subset E$ with $rk(F)=s < r$ such that
\begin{eqnarray*}
\frac{\chi(E(n))}{r}-\frac{\chi(F(n))}{s}
&=& \left(\frac{c_1(E)\cdot H}{r}-\frac{c_1(F)\cdot H}{s}\right)n\\
&&+\left(\frac{ch(E)\cdot td(X)}{r}-\frac{ch(F)\cdot td(X)}{s}\right)\\
&=&0
\end{eqnarray*}
for all $n>>0$.
This  contradicts the gcd condition.
\end{proof}

Suppose $X$ is an abelian (respectively K3) surface such that
  $ch$ (respectively $ch(1+[X])$) is primitive.
Then $M_X^s(r,c_1, c_2, H)=M_X(r,c_1, c_2, H)$ provided that
$H$  sufficiently general (this can be deduced from \cite{yoshi}).

\begin{thm}\label{thm:modsurfaces}
Let $X, H\ldots$ be as above, and let $M=M_X(r,c_1, c_2, H)$ and
$M^s=M_X^s(r,c_1, c_2, H)$. Then
\begin{enumerate}
\item The  generalized Hodge conjecture holds for every component
of  $M$ if $X$ is rational,
$-K\cdot H<0$, and the hypothesis of lemma~\ref{lemma:gcdofrc1chi}
holds.

\item The  Hodge (respectively generalized Hodge) conjecture holds for every component of
$M$ if  $-K\cdot H<0$, the hypothesis of lemma~\ref{lemma:gcdofrc1chi}
 holds and $X$ is  minimal ruled over a curve $C$ all of whose powers
satisfy the  Hodge (respectively generalized Hodge) conjecture.

\item The  Hodge (respectively generalized Hodge) conjecture holds for every component of
$M$ if  $M=M^s$,  $X$ is abelian or K3, and  all powers of $X$ satisfy the 
 Hodge (respectively generalized Hodge) conjecture.
\end{enumerate}
\end{thm}

\begin{proof}
The result will be reduced to proposition~\ref{prop:ABtype} in 
each of the above cases. Note that
 $X^k$ satisfies the (generalized) Hodge conjecture for all $k$ either by hypothesis
or by proposition \ref{prop:hodgesurface}.
Under the assumptions of (1) or (2),  Maruyama \cite{maruyama} has
shown that $M$ that there
is a universal sheaf $\E$ on $X\times M$ and that 
$Ext_X^2(E_1,E_2) = 0$ for any two sheaves in $E_i\in M$.
Therefore, Beauville's theorem \cite{beauville}
applies to show that
$H^*(M)$ is generated by the K\"unneth components of $c_i(\E)$
as an algebra. The result follows by proposition~\ref{prop:ABtype}.
This finishes the proof in cases (1) and (2).

In case (3), we can apply the main theorem of Markman \cite{markman}
to see that $H^*(M)$ is generated by the K\"unneth components of Chern
classes of  a quasi-universal sheaf $\E$ on $X\times M$.
\end{proof}

\begin{cor}
The Hodge conjecture holds for $M$ when $X$ is an abelian surface
or a K3 surface satisfying the conditions of
proposition~\ref{prop:hodgesurface} (e)   or (f).
\end{cor}

\begin{proof}
Follows from proposition~\ref{prop:hodgesurface}.
\end{proof}

\begin{cor}
    The generalized Hodge conjecture holds for $M$ when $X$ is
 a product of two elliptic curves without complex multiplication
 or $X$ is a simple abelian surface of type I or II.
\end{cor}

\begin{proof}
    This follows from corollary~\ref{cor:Hazama}
 \end{proof}

\section{Arithmetic analogues}

Let $k$ be a field which is finitely generated over a prime field.
Let $\bar k$ denote the separable closure  and $G = Gal(\bar k/k)$
the  Galois group. Choose a prime $l\not= char\, k$.
Given a variety $X$ defined over $k$, let
 $\bar X = X\times_{spec k} spec\, \bar k$.
Then the \'etale cohomology group
 $H_{et}^i(\bar X, \Q_l)$
is a $\Q_l$-vector space with a continuous $G$-action.
Tate twisting in this context amounts twisting by a
power of the cyclotomic character (see \cite{katz} for
a rapid introduction to these ideas). 
When $X$ is smooth and projective we refer to an invariant in
$H_{et}^{2i}(\bar X, \Q_l(i))^G$ as a Tate cycle.
 The fundamental class of a subvariety
defined over $k$ is a Tate cycle. The Tate conjecture \cite{tate1, tate3}
 claims  conversely that
 the  space of Tate cycles is  spanned by these 
fundamental classes. This can be viewed as an analogue of
the Hodge conjecture. To make this analogy clearer, note
that the space of Hodge cycles on the $2i$th cohomology of
a complex smooth projective variety  is isomorphic to  
$Hom_{MHS}(\Q, H^{2i}(X,\Q)(i))$. Similarly the
space of Tate cycles is $Hom_{G}(\Q_l, H_{et}^{2i}(\bar X,\Q_l(i)))$.

 Deligne \cite[I]{dm} has proposed a variant of the Hodge conjecture 
that says roughly that the property of being a Hodge cycle should be invariant
under the automorphism group of $\C$.
We will give the precise formulation in a manner which generalizes
easily. Suppose that $k$ is a finitely generated field of
characteristic $0$.
Jannsen \cite[\S2]{jannsen} has constructed the abelian category
$MR_k$ of  mixed realizations. An object in this category 
is a collection of the following data:
\begin{enumerate}
\item A bifiltered finite dimensional $k$-vector space $(H_{dR}, F,W)$
\item A filtered  finite dimensional $\Q_l$ vector space
with a continuous $G$-action $(H_l, W)$ for
  each prime $l$.
\item A $\Q$-mixed Hodge structure $H_\sigma$, $(H_\sigma,F,W)$
for each  embedding $\sigma:\bar k\to \C$.

\item Comparison isomorphisms $H_\sigma\otimes \Q_l\cong H_l$
and $H_{dR}\otimes\C \cong H_\sigma\otimes \C$ respecting the
filtrations.
\end{enumerate}
The morphisms are constructed so that the obvious projections
of $MR_k$ are functors.
In particular, for each
embedding $\bar\sigma:\bar k\to \C$, there is a functor from
$\Phi_{\bar\sigma}:MR_k\to MHS$.
For example, if $X$ is a smooth projective variety over $k$, the
collection of de Rham, $l$-adic and analytic cohomologies
$$(H^i_{dR}(X), H_{et}^i(\bar X,\Q_l), H^i(\bar X\times_{ \sigma}
spec\, \C,\Q)),$$
 along  with the comparison isomorphisms
 and appropriate filtrations provides an example
of an object $H^{i}_{AH}(X,0)\in MR_k$. More generally, 
Jannsen \cite[6.11.1]{jannsen}
has a constructed a (co)homology theory
$H^*_{AH}(X,j),\Hbm^{AH}_*(X,j)$ from the category of $k$-varieties
$X$ to $MR_k$ such that 
$$\Phi_\sigma (H^*_{AH}(X,j)) = H^*(\bar X\times_\sigma  spec\,\C,\Q(j))$$
$$\Phi_\sigma (\Hbm^{AH}_*(X,j)) = \Hbm_*(\bar X\times_\sigma  spec\,\C,\Q(j)).$$
Given an object $H\in MR_k$, the space of absolute Hodge cycles
$$\Gamma(H) =
\{(\alpha_{dR}, \ldots)\mid
\mbox{the components are compatible, }
\alpha_{dR}\in F^0H_{dR}\cap W_0H_{dR}\} $$
If $X$ is a smooth projective variety over $k$,
then each component $\alpha_\sigma$ of $\alpha \in \Gamma(H^{2i}_{AH}(X,i))$
is a Hodge cycle. Deligne's conjecture, which we will
refer to  as the absoluteness conjecture, is that for each $\sigma$,
any  Hodge cycle $H^{2i}(\bar X\times_\sigma  spec\,\C,\Q(i))$ arises
as $\alpha_\sigma$ for some absolute Hodge cycles $\alpha$.
This conjecture would follow from the Hodge conjecture since
the collection of fundamental classes of an algebraic cycle defined
over  $k$  for the above cohomology theories is  an absolute
Hodge cycle.

We need the analogues of the results of section 1.

\begin{lemma}\label{lemma:lifttate}
Let  $X$ and $Y$ be smooth projective $k$-varieties 
and suppose that 
$f:H_{et}^{2i}(\bar X,\Q_l(i))\to H_{et}^{2j}(\bar Y,\Q_l((j))$
 is a surjective
morphism  induced by a correspondence defined over $k$.
If the Tate conjecture holds for $X$ and if 
$char\, k = 0$  or  $H^*(\bar X)$ is a semisimple $G$-module, then the
Tate conjecture  holds for $Y$. If $char\, k = 0$ and
the absoluteness conjecture holds for
$X$, then it holds for $Y$.
\end{lemma}

\begin{proof} If the map $f:H^*(\bar X)\to H^*(\bar Y)$ of Galois modules
splits, then any Tate cycle $\alpha$ in $H^*(\bar Y)$ can be lifted to a
Tate cycle $\beta$ in $H^*(\bar X)$.   $\beta$ would be algebraic
if Tate's conjecture held for $X$, therefore $\alpha$ is also algebraic.
 The splitting of $f$ is
immediate when $H^*(\bar X)$ is semisimple. When   $char \, k = 0$, 
the splitting follows from  \cite[1.2]{jannsen}.

Similarly by corollary \ref{cor:lifthodge}, a Hodge cycle 
 $\alpha$ in  $H^*(\bar Y\times_\sigma spec\,\C)$ can be lifted to a Hodge cycle 
$\beta$ in $H^*(\bar X\times_\sigma spec\,\C)$. Assuming the absoluteness
 conjecture for  $X$, $\beta$
and hence $\alpha$,  would necessarily extend an absolute Hodge cycle.
\end{proof}

\begin{remark}
It is a conjectured that $H^*(\bar X)$ is always semisimple. However
this has been established in only very special cases
\cite{faltings, tate1.5, zarhin}.
\end{remark}

Jannsen \cite[7.3]{jannsen} has also formulated a version of Tate's
conjecture for singular quasiprojective varieties using Borel-Moore
\'etale homology which can be defined as dual to compactly supported
\'etale cohomology as above.
Given such a variety $X$, the fundamental class
of any $i$ dimensional $k$-subvariety lies in $\Hbm^{et}_{2i}(\bar X,\Q_l(-i))^G$.
We say that the Tate conjecture holds for $X$ if  $\Hbm^{et}_{2*}(\bar
X,\Q_l(*))^G$
is spanned by these classes. The absoluteness conjecture can
 be extended in a similar fashion: namely, that Hodge cycles in 
$\Hbm_{2i}(\bar X\times_\sigma  spec\,\C,\Q(-i))$ lift to absolute Hodge cycles
$\Gamma(\Hbm^{AH}_*(X,-i))$.

\begin{prop}[Jannsen\cite{jannsen}]\label{prop:jannsentate}
If $char\, k = 0$, the Tate conjecture holds for a $k$-variety
$U$ if it holds for a desingularization of a compactification of $U$.
\end{prop}

\begin{lemma}\label{lemma:jannsenAC}
Suppose that $k$ is an algebraically closed field of characteristic
 $0$. The absoluteness conjecture holds for a $k$-variety
$U$ if it holds for a desingularization of a compactification of $U$.
\end{lemma}

\begin{proof}
Fix an embedding $\sigma:k\to \C$, and let $Y_\sigma = \bar
Y\times_\sigma spec\, \C$.
Let $\tilde X$ be a  desingularization of a compactification of $U$
satisfying the absoluteness conjecture. 
Thus the map from the space of absolute Hodge to   Hodge cycles
$$\Gamma H_{2i}^{AH}(\tilde X, i)\to  H_{2i}^{Hodge}(\tilde X_\sigma)$$
is surjective. 
Let $X \subseteq \tilde X$ be the preimage of $U$.
Then the composition of 
$$H_{2i}(\tilde X_\sigma,\Q(i))\to H_{2i}( X_\sigma,\Q(i))\to 
H_{2i}(U_\sigma,\Q(i))$$
induces a surjection
$$H_{2i}^{Hodge}(\tilde X_\sigma,\Q(i))\to H_{2i}^{Hodge}(U_\sigma,\Q(i))$$
(see \cite[pp 113-114]{jannsen}).
Consequently 
$$\Gamma H_{2i}^{AH}(U, i)\to  H_{2i}^{Hodge}(U_\sigma, i)$$
is surjective since it factors through the composition
$$\Gamma H_{2i}^{AH}(\tilde X,\Q(i))\to H_{2i}^{Hodge}(U_\sigma,\Q(i)).$$
\end{proof}

\begin{lemma}\label{lemma:strattate}
Suppose that $char\, k = 0$.
Let $Z\subset X$ be a closed subset of a projective variety
$X$, and let $U = X-Z$. The Tate conjecture (respectively
 the absoluteness conjecture) holds for
$X$ if it holds for $U$ and $Z$.
\end{lemma}

\begin{proof}
There is an exact sequence of Galois modules
$$\Hbm_{2k}^{et}(\bar Z,\Q_l)\to \Hbm_{2k}^{et}(\bar X,\Q_l)\to \Hbm_{2k}^{et}(\bar U,\Q_l) .$$
From weight considerations, the image $I$ of
$\Hbm_{2k}^{et}(\bar Z,\Q_l)$ in $\Hbm_{2k}^{et}(\bar X,\Q_l)$ coincides with
the image of $\Hbm_{2k}^{et}(\bar{\tilde  Z},\Q_l)$ for any desingularization
$\tilde Z\to Z$.
By \cite[1.2]{jannsen},  $I$ possesses a complement in the category
of  Galois modules.
It follows that there is an exact sequence of Tate cycles
$$\Hbm_{2k}^{et}(\bar Z,\Q_l(-k))^G\to \Hbm_{2k}^{et}(\bar X,\Q_l(-k))^G\to 
\Hbm_{2k}^{et}(\bar U,\Q_l(-k))^G $$

Likewise by lemma \ref{lemma:strat}, there is an exact sequence
of absolute Hodge cycles 
$$\Gamma(\Hbm_{2k}(Z,-k))\to \Gamma(\Hbm_{2k}(X,-k))\to 
\Gamma(\Hbm_{2k}(U, -k)). $$
The rest of the argument proceeds exactly as in the proof of
lemma \ref{lemma:strat}.
\end{proof}

From this lemma, we can deduce the analogue of
corollary \ref{cor:strat}  for the Tate
and Deligne conjectures when
$k$ has characteristic $0$. 

\begin{lemma}\label{lemma:fiber3}
Suppose $char\, k =0$.
 Let $f:X\to Y$ be  a morphism of smooth $k$-varieties
which is a Zariski locally trivial fiber
 bundle with fiber $F$. Suppose that $F$ is  smooth
and  that $\Hbm_*^{et}(\bar F, \Q_l)$ is spanned by algebraic
cycles. Then  the Tate (or absoluteness)  conjecture
holds for $X$ if it holds for $Y$.
\end{lemma}

The proof of this is virtually identical to the proof of lemma 
\ref{lemma:fiber}
(the Leray-Hirsch theorem for etale cohomology can be deduced
by the comparison theorem).

 Milne \cite[sect. 4]{milne} has proven a version of 
 theorem \ref{thm:HdgLef} for the Tate conjecture. We note a special case.

\begin{prop} If $X$ is a polarized abelian variety over $k$ such
that the image of $G$ is Zariski dense in $GSp(H_{et}^1(\bar X,\Q_l))$,
then the Tate conjecture holds for all powers of $X$.
\end{prop}

\begin{proof} In the terminology of \cite{milne}, the Zariski
closure of the image of
$G$ is a subgroup of the Tate group which is contained in the
Lefschetz group $\subseteq GSp(H_{et}^1(\bar X,\Q_l))$.
When equality holds, Tate's conjecture holds for $X$ and its
powers.
\end{proof}

\begin{cor}\label{cor:endZ}
If $k$ is a number field, $End(\bar X) = \Z$ and
$dim\, X$ is either odd or equal to $2$
or $6$, then  the Tate conjecture holds for all powers of $X$.
\end{cor}

\begin{proof}
Under these conditions $G$ is dense in $GSp(H^1(X))$ by
 a theorem of Serre \cite[2.28]{serre} (see also \cite[6.1]{chi}).
\end{proof}

Deligne \cite[I,II]{dm} has proven that the absoluteness conjecture holds
 for certain special classes of varieties.

\begin{thm}[Deligne] Let $X$ be an abelian variety, a product of
 K3 surfaces, or a product of Fermat hypersurfaces. Then
the absoluteness conjecture holds.
\end{thm}

\begin{cor}\label{cor:H2AHcurve}
Let $X$ be smooth projective curve, then 
the absoluteness conjecture holds for all powers of $X$.
\end{cor}

\begin{proof}
The argument is similar to the  proof of proposition \ref{prop:jacobian}.
\end{proof}

 The following are analogues of proposition~\ref{prop:ABtype}
and theorem \ref{thm:powersXtoM}.

\begin{prop}\label{prop:ABtypeTate}
    Suppose that $Y$ and $N$ are smooth projective varieties  over
 $k$ such that there exists a finite collection of algebraic correspondences 
 on $Y\times N$ such their K\"unneth components generate
 the cohomology ring $H_{et}^*(\bar N,\Q_l)$.  Assume that
either $H_{et}^1(\bar Y,\Q_l)$ is semisimple or $char\, k = 0$.
If such that $Y^m$ satisfies the  Tate conjecture
for all $m \le dim N$, then the  Tate  conjecture holds for $N$.
If $k$ is an algebraically closed field of characteristic $0$,
then the absoluteness conjecture for  $Y^m$ for all $m\le dim N$
implies the absoluteness conjecture for $N$.
 \end{prop} 

 \begin{proof}
   As in the proof of proposition~\ref{prop:ABtype} we get a
   surjection
$$A: \oplus_{J,I} H^{*}(\bar Y^*)(*) \to H^a(\bar N)$$
induced by a correspondence,
where $H^*$ is taken to be $l$-adic cohomology or $H^*_{AH}$.
An appeal to lemma~\ref{lemma:lifttate} finishes the proof.
 \end{proof}

\begin{thm}
Let $X$ be a smooth projective curve over $k$, and
let $M=U_X(n,d)$ with $n$ and 
$d$ coprime. Assume that
either $H_{et}^1(X,\Q_l)$ is semisimple or $char\, k = 0$.
Then the Tate conjecture holds for $M$ if it holds for all powers
of $X$.
If $k$ is an algebraically closed field of characteristic $0$,
then the absoluteness conjecture holds for $M$.
\end{thm}

\begin{remark} The part of the theorem concerning the Tate conjecture
in characteristic $0$ is due to del Ba\~no \cite{delb}.
\end{remark}

\begin{proof}
The argument is identical to second proof of theorem
\ref{thm:powersXtoM} with proposition~\ref{prop:ABtypeTate}
replacing proposition~\ref{prop:ABtypeTate}. (Note that
Beauville's proof \cite{beauville} is algebraic and is
valid in positive characteristic provided etale cohomology
is used in place of singular cohomology.)
\end{proof}

\begin{cor}
The Tate conjecture holds for $M$ if the hypotheses of corollary
\ref{cor:endZ} are satisfied for $k$ and $J(X)$.
\end{cor}

\begin{proof}
As in the proof of proposition \ref{prop:jacobian}, one gets a surjection
$H_{et}^*(\bar J(X)^m)\to H^*_{et}(\bar X^m)$.
Lemma \ref{lemma:lifttate} implies the Tate conjecture for $X^m$.
\end{proof}

The proof of theorem \ref{thm:hilb} can be modified to yield:

\begin{thm} Suppose that $char\, k = 0$.
Let $X$ be a smooth projective surface such that all powers $X^m$,
with $m\le n$, satisfy the Tate conjecture 
(respectively the absoluteness conjecture).
Then $X^{[n]}$ satisfies the Tate conjecture (respectively
the absoluteness conjecture).
\end{thm}

\begin{cor} Suppose that $char\, k =0$. Let $X$ be a smooth projective surface
over $k$ with Kodaira dimension $\kappa(X)\le 0$,  then
 the absoluteness conjecture holds for $X^{[n]}$ for any $n$.
\end{cor}

\begin{proof}
First note that when $X$ is a product of $\PP^1$ and a
smooth projective curve, the absoluteness conjecture holds for
all powers of $X$ by corollary \ref{cor:H2AHcurve}
and  lemma \ref{lemma:fiber3}.
Then arguing as in the proof of \ref{prop:hodgesurface}, we see that
if there is   a smooth projective surface
$Y$, all of whose powers satisfy the absoluteness conjecture,
and a dominant rational map $Y\dashrightarrow X$, then
the absoluteness conjecture holds for all powers of $X$ as well.

If $X$ is a smooth surface of Kodaira dimension zero,
it is a:
\begin{enumerate}
\item rational or  ruled,
\item abelian or bielliptic, or
\item  K3 or Enriques 
\end{enumerate}
surface \cite{beauville-surface}.
Then there exists a dominant map $Y\dashrightarrow X$ where
$Y$ is a product of $\PP^1$ and curve in case 1,
an abelian surface in case 2, or a K3 surface in case 3.
This finishes the proof.
\end{proof}

We have an analogue of theorem~\ref{thm:modsurfaces}. For simplicity,
we state the most interesting part.

\begin{thm}
   Suppose that $char\, k = 0$.
Let $X$ be an abelian or K3 surface defined over $k$,
and assume that classes $r,c_1, c_2, H$ have been chosen so that
 $M=M_X(r,c_1, c_2, H)=M_X^s(r,c_1, c_2, H)$.
If the Tate conjecture holds for all powers of $X$, then it holds for
$M$. If $k$ is algebraically closed, then the absoluteness conjecture holds for $M$.
\end{thm}

\begin{proof}
  The proof is identical to the proof of theorem~\ref{thm:modsurfaces}
with proposition~\ref{prop:ABtypeTate} in the place of  
proposition~\ref{prop:ABtype}.
\end{proof}


\begin{thebibliography}{AB}


\bibitem[AO]{AO} D. Abramovich, F. Oort, 
{\em  Alterations and resolution of singularities.}
 Resolution of singularities (Obergurgl, 1997),
 39--108, Progr. Math., 181, Birkh"user, (2000)

\bibitem[AB]{AB} M. Atiyah, R. Bott, {\em Yang-Mills equations
over Riemann surfaces}, Phil. Trans. Royal Soc. London 308,  523-615(1983)

\bibitem[Be1]{beauville-surface} A. Beauville, {\em
Surface alg\'ebriques complexes}, Asterisque (1978)

\bibitem[Be2]{beauville} A. Beauville, {\em  Sur la cohomologie de
 certains espaces de modules de fibr\'es vectoriels},
Geometry and analysis, 37--40, Tata Inst. (1995)

\bibitem[BB]{bb} A. Bialynicki-Birula, {\em Some theorems
on actions of algebraic groups}, Ann. of Math. (1972)


\bibitem[BGL]{bgl} E. Biffet, F. Ghione, M. Letizia,
{\em On the Abel-Jacobi map for divisors of higher rank on
a curve}, Math. Ann. 299 (1994)


\bibitem[BN]{BN} I. Biswas, M. S. Narasimhan, {\em Hodge classes of 
moduli spaces of parabolic bundles over the general  curve.},
J. Alg. Geom. 4, 697--715. (1997)

\bibitem[CM1]{cm1} M. de Cataldo, L. Milgiorini,
{\em The Douday space of a complex surface}, Advances in Math (2000)

\bibitem[CM2]{cm2} M. de Cataldo, L. Milgiorini,
{\em Chow groups and the motive of the Hilbert scheme of points on
a surface}, J, Algebra (to appear)


\bibitem[C]{chi} W. Chi, {\em $l$-adic and $\lambda$-adic 
representations associated to Abelian varieties over number fields}
Amer. J. Math 114, 315-353 (1992)

\bibitem[CN]{conte-murre} A. Conte, J. Murre, {\em The Hodge
conjecture for fourfolds admitting a covering by rational curves},
Math. Ann. 238 (1978)

\bibitem[db]{delb} S. del Ba\~no, {\em Chow motive of some moduli
    spaces}, Crelles J. 532 105-132 (2001)



\bibitem[D1]{deligne-hodge} P. Deligne, {\em Th\'eorie de Hodge II},
  Pub. IHES 40. 5--57 (1971)

\bibitem[D2]{deligne} P. Deligne, {\em La conjecture de Weil pour les
surfaces K3}, Invent. Math. 15  206--226 (1972)

\bibitem[D3]{deligne2} P. Deligne {\em La conjecture de Weil I}
  Publ. IHES 43, (1974), 273--307.
 
\bibitem[DMOS]{dm} P. Deligne, J. Milne, A. Ogus, K. Shi,
 {\em Hodge cycles, motives and Shimura varieties}, LNM 900, Springer-Verlag
(1982)

\bibitem[Fa]{faltings} G. Faltings, {\em Endlichkeitss\"atze f\"ur
abelsche  Variet\"aten \"uber Zahlk\"orpern. } Invent. Math. 73, 349--366.(1983)

\bibitem[F1]{fulton} W. Fulton, {\em Intersection theory},
Springer-Verlag (1984)



\bibitem[FH]{FH} W. Fulton, J. Harris, {\em Representation theory}
Springer-Verlag (1991)


\bibitem[Gi]{gieseker} D. Gieseker, {\em On the modulis of vector
bundles on an algebraic surface}, Ann. Math. 106 (1977)

\bibitem[Grd]{gordon} B. Gordon, {\em The Hodge conjecture for
abelian varieties}, Appendix to Survey of the Hodge Conjecture 2nd
ed. by J. Lewis, AMS (1999)

\bibitem[GM]{gm} M. Goresky, R. Macpherson, {\em Stratified Morse theory}
Springer-Verlag (1980)

\bibitem[G]{gott} L. G\"ottsche, {Hilbert schemes of Zero dimensional
subschemes of smooth varities}, Lect. Notes Math 1572, 
Springer-Verlag (1994)

\bibitem[Gr]{groth} A. Grothendieck, Techniques de construction et théorèmes d'existence en géométrie algébrique. IV., Sem. Bourbaki, Exp 221,
(1960)

\bibitem[Gr2]{groth2} A. Grothendieck, Hodge's general conjecture
is false for trivial reasons, Topology 8, (1969)

\bibitem[Ha]{hain} R. Hain, {\em  Moduli of Riemann surfaces, transcendental aspects.} School on Algebraic Geometry (Trieste,
1999), 293--353, ICTP Lect. Notes, 1, Abdus Salam Int. Cent. Theoret. Phys., Trieste, (2000)

\bibitem[Hz]{hazama} F. Hazama, {\em The generalized Hodge conjecture
for stably nondegenerate abelian varieties}, Compositio Math. 93,
(1993)


\bibitem[H]{hodge} W. Hodge, {\em The topological invariants
of algebraic varieties} , Proc. ICM (1950)

\bibitem[J]{jannsen} U. Jannsen, {\em Mixed motives and
algebraic K-theory}, Lect notes in math 1400, Springer-Verlag (1990)

\bibitem[K]{katz} N. Katz, {\em Review of $l$-adic cohomology}
 Motives (Seattle, WA, 1991), Proc.
Sympos. Pure Math., 55, Part 1, Amer. Math. Soc., 21-30 (1994)

\bibitem[L]{lewis} J. Lewis, A survey of the Hodge conjecture,
CRM Monographs 10, AMS (1999)

\bibitem[K]{sga7} N. Katz, {\em Etude cohomologie des pinceaux de
 Lefschetz}  SGA7, exp. XVIII, LNM  340, Springer-Verlag (1973)

\bibitem[Mrk]{markman} E. Markman, {\em Generators of the cohomology 
ring of moduli spaces of sheaves on symplectic surfaces}


\bibitem[Ma]{maruyama} M. Maruyama, {\em Moduli of stable
sheaves, II}, J. Math. Kyoto 18 (1978)

\bibitem[Mi]{milne} J. Milne, {\em Lefschetz classes on abelian
varieties} Duke Math. J. (1999)

\bibitem[MZ]{moonzar} B. Moonen, Y. Zarhin, {\em Hodge classes on 
abelian varieties of low dimension}
Math. Ann (1999)

\bibitem[Mk1]{mukai} S. Mukai, {\em Symplectic structure of the moduli
    space of 
sheaves on an abelian or $K3$ surface} Invent. Math.  77  (1984),
    101-116

\bibitem[Mk2]{mukai2} Mukai, S. {\em On the moduli space of bundles on $K3$
 surfaces. I}  Vector bundles on algebraic varieties (Bombay, 1984), 
 341--413, Tata Inst. ( 1987)

\bibitem[Mu]{mumford} D. Mumford, {\em Abelian varieties}
Tata Inst. 


\bibitem[Mu1]{murty2} K. Murty, {\em Computing the Hodge
group of an abelian variety}, S\'em. Th\'eorie Nom. Paris 1988-89,
Birkh\"auser (1990)

\bibitem[Mu2]{murty} K. Murty, {\em Hodge and Weil classes on
Abelian varieties}, in The arithmetic and geometry of algebraic
cycles, Kluwer (2000)



\bibitem[R]{ribet} K. Ribet, {\em Hodge classes on certain
types of abelian varieties}, Amer. J. Math. 105 (1983)


\bibitem[Sc]{schoen} C. Schoen. {\em Varieties dominated by products}
Int. J. Math 7, 541--571 (1996)

\bibitem[Se]{serre} J. P. Serre, {\em Resum\'e de cours 1984-1985},
Ouevres IV, Springer-Verlag (2000)

\bibitem[Sh]{shioda} T. Shioda, {\em Algebraic cycles on
Abelian varieties of Fermat type} Math. Ann. 258 (1981)

\bibitem[Sp]{spanier} E. Spanier, {\em Algebraic Topology}
McGraw-Hill (1966)

\bibitem[T1]{tate1} J. Tate, {\em Algebraic cycles and poles of zeta
functions} , Arithmetic Algebraic Geometry 93--110, Harper and Row (1965)

\bibitem[T2]{tate1.5} J. Tate, {\em Endomorphisms of abelian varieties
 over  finite fields. } Invent. Math. 2,  134--144.(1966)

\bibitem[T3]{tate3}J. Tate,{\em Conjectures
on algebraic cycles in $l$-adic cohomology}
 Motives (Seattle, WA, 1991), Proc.
Sympos. Pure Math., 55, Part 1, Amer. Math. Soc.,71-83,  (1994)

\bibitem[V]{verdier} J. L. Verdier, {\em Ind\'ependance par rapport
\`a $l$ de polyn\^omes charact\'eristiques ...}, Sem. Bourbaki 423, LNM 383, Springer-Verlag (1974)

\bibitem[Y]{yoshi} Yoshioka, K. {\em Chamber structure of
    polarizations and the moduli of stable sheaves on a ruled
    surface.} Internat. J. Math.  7  (1996)

\bibitem[Z1]{zarhin}  Y. Zarhin,
{\em Endomorphisms of Abelian varieties over fields of finite
  characteristic.} Izv. Akad. Nauk SSSR Ser. Mat. 39 (1975)

\bibitem[Z2]{zarhin2}  Y. Zarhin,
{\em Hodge groups of $K3$ surfaces}
 J. Reine Angew. Math.  341, 193-220, (1983)
  
\end{thebibliography}
\end{document}